\documentclass[11pt]{amsart}
\usepackage{amssymb,amsmath,amsthm,newlfont}
\newtheorem{prop}{Proposition}[section]
\newtheorem{thm}[prop]{Theorem}
\newtheorem{lem}[prop]{Lemma}
\newtheorem{cor}[prop]{Corollary}
\newtheorem{definition}[prop]{Definition}

\theoremstyle{remark}
\newtheorem{rem}[prop]{Remark}
\newtheorem{ex}[prop]{Example}

\begin{document}

\title[Amenability of $\cal A(X)$]{Amenability of algebras
\\ of approximable operators}

\author{A.~Blanco}
\address{Department of Pure Mathematics,
Queen's University Belfast, Belfast BT7 1NN, UK}
\email{a.blanco@qub.ac.uk}

\author{N.~Gr{\o}nb{\ae}k}
\address{Institute for Mathematical Sciences, Department of
Mathematics, Universitetsparken~5, DK-2100 Copenhagen {\O},
Denmark}
\email{gronbaek@math.ku.dk}




\begin{abstract}
We give a necessary and sufficient condition for amenability of
the Banach algebra of approximable operators on a Banach space. We
further investigate the relationship between amenability of this
algebra and factorization of operators, strengthening known
results and developing new techniques to determine whether or not
a given Banach space carries an amenable algebra of approximable
operators. Using these techniques, we are able to show, among
other things, the non-amenability of the algebra of approximable
operators on Tsirelson's space.
\end{abstract}

\maketitle

\section{Introduction}

Let $\cal A$ be a Banach algebra and let $\cal X$ be a Banach
space which is also an $\cal A$-bimodule. Then $\cal X$ is a
Banach $\cal A$-bimodule if there exists a constant $M$ so that
$\|a\cdot x\|\leq M\|a\|\|x\|$ and $\|x\cdot a\|\leq M\|a\|\|x\|$
$(a\in\cal A,\,x\in\cal X)$. A (continuous) derivation from $\cal
A$ to $\cal X$ is a (bounded) linear map $D:\cal A\to\cal X$ that
satisfies the identity
\[
D(ab) = D(a)\cdot b + a\cdot D(b)\qquad (a,b\in\cal A).
\]
Every map of the form $a\mapsto a\cdot x - x\cdot a$ $(a\in\cal
A)$, where $x\in\cal X$ is fixed, is a continuous derivation.
Derivations of this form are called inner derivations.

If $\cal X$ is a Banach $\cal A$-bimodule, then its topological
dual, $\cal X^*$, is also a Banach $\cal A$-bimodule under the
actions
\[
(a\cdot f)(x) := f(x\cdot a)\quad \text{ and }\quad (f\cdot a)(x)
:= f(a\cdot x)\qquad (a\in\cal A,\,x\in\cal X,\,f\in\cal X^*).
\]

The Banach algebra $\cal A$ is said to be \emph{amenable} if, for
every Banach $\cal A$-bimodule $\cal X$, every continuous
derivation $D:\cal A\to\cal X^*$ is inner.

For example, the group algebra, $L^1(G)$, of a locally compact
group is amenable if and only if the group $G$ is amenable
\cite{J1}; a $C^*$-algebra is amenable if and only if it is
nuclear \cite{Co,Ha}; and a uniform algebra on a compact Hausdorff
space $\Omega$ is amenable if and only if it is $C(\Omega)$
\cite{S}.

In this note we shall be concerned with the amenability of the
algebra $\cal A(X)$ of \emph{approximable operators} on a Banach
space $X$, i.e., the operator norm closure in $\cal B(X)$ of the
ideal $\cal F(X)$ of continuous finite-rank operators on $X$,
where $\cal B(X)$ denotes the algebra of all bounded linear
operators on $X$. (When $X$ has the approximation property, $\cal
A(X)$ coincides with the ideal of compact operators on $X$.) In
this setting the main problem is to characterize amenability of
$\cal A(X)$ in terms of properties of $X$.

The study of amenability of $\cal A(X)$ goes back to \cite{J1},
where it is shown that $\cal A(X)$ is amenable if $X = \ell_p$ for
$p\in(1,\infty)$, or $X = C[0,1]$. Further progress in the study
of amenability of this algebra is made in \cite{GJW}. In this last
paper a geometric property, called property $\Bbb A$, is
introduced, and it is shown that Banach spaces with this property
carry amenable algebras of approximable operators. Banach spaces
with property $\Bbb A$ include all classical Banach spaces, $\cal
L_p$-spaces $(1\leq p\leq\infty)$, spaces with a subsymmetric,
shrinking basis, and certain kinds of tensor product of Banach
spaces with property $\Bbb A$.

In this note we continue the study of amenability of the algebra
$\cal A(X)$. Building upon ideas from \cite{GJW} we shall develop
new techniques that will allow us not only to improve several
results from \cite{GJW} but also to answer some of the questions
left open there. In particular, we will show that the algebra of
approximable operators on Tsirelson's space is not amenable. An
important fact that should become apparent throughout these pages
is that a full understanding of amenability of $\cal A(X)$ will
necessarily rely on a good understanding of the finite-dimensional
case.

The paper has been organized as follows. In the next section, we
have gathered some terminology and basic facts we need. In
Section~3, we give a necessary and sufficient condition for
amenability of $\cal A(X)$. Unfortunately, practical use of this
condition depends on our capability to find good estimates for the
projective norm of certain elements called generalized diagonals.
In Section~4, we follow a different approach. The results of this
section are to a great extent motivated by the notion of
approximate primariness introduced in \cite{GJW}. We explore some
of the ideas behind this notion, specially, its connection with
factorization properties of operators. Finally, in Section~5, we
establish the non-amenability of $\cal A(X)$ for every Banach
space $X$ in a certain family of Tsirelson-like~spaces. In doing
this we shall rely on results from the previous sections.

\section{Preliminaries}

In this section we have gathered some notations and basic results
that we shall use throughout these pages.

To simplify the statement of the results, we shall denote by
$\ell_\infty$ the linear space, usually denoted by $c_0$, of all
bounded scalar sequences tending to zero. Given a normed space $X$
we denote by $X^*$ its topological dual. If $X$ and $Y$ are
isomorphic (resp.~isometric) normed spaces, we write this as
$X\simeq Y$ (resp.~$X\cong Y$), and denote by $d(X,Y)$ the
Banach--Mazur distance between them, that is, the infimum of
numbers $\|T\|\|T^{-1}\|$, where $T$ is an isomorphism between $X$
and $Y$.

The adjoint of an operator $T:X\to Y$ is denoted by $T^*$ and we
write $\text{rg}\,T$ (resp.~$\text{rk}\,T$) for the range
(resp.~rank) of $T$. The identity operator on a normed space $X$
is denoted by $I_X$ or just $I$ if the space $X$ is clear from
context.

By the inversion constant of a surjective linear map, $Q:X\to Y$,
between Banach spaces we mean the operator norm of the inverse of
the linear isomorphism $\widetilde{Q}:X/\ker{Q}\to Y$ induced by
$Q$, that is, $\|\widetilde{Q}^{-1}\|$. Given Banach spaces $X$,
$Y$ and $Z$, a bilinear map $\varphi:X\times Y\to Z$ will be said
to be $M$-open if for every $z\in Z$ there exist $x\in X$ and
$y\in Y$ such that $\varphi(x,y) = z$ and $\|x\|\|y\|\leq M\|z\|$.

We write $\|\cdot\|_\wedge$ (resp. $\|\cdot\|$) for the projective
(resp. operator) norm. If two norms, $\|\cdot\|_1$ and
$\|\cdot\|_2$, on a linear space are equivalent we write this as
$\|\cdot\|_1 \sim \|\cdot\|_2$. Given a set of vectors $\{e_i :
i\in I\}$ in a Banach space, we denote by $[e_i]_{i\in I}$ the
closure of its linear span.

Let $\mathbf{e} = (e_i)$ be a 1-unconditional basis for the Banach
space $(E,\|\cdot\|)$, and let $(X_i,\|\cdot\|_i)$ be a sequence
of  Banach spaces. We let
\[
\Big(\bigoplus\nolimits_i X_i\Big)_{\mathbf{e}} =
\Big\{(x_i)\in\prod\nolimits_i X_i : \sum\nolimits_i \|x_i\|_i e_i
\; \text{converges in $E$}\Big\},
\]
endowed with the norm $\|(x_i)\| := \big\|\sum_i \|x_i\|_i
e_i\big\|$. It is well known that $(\bigoplus_i X_i)_{\mathbf{e}}$
is a Banach space. Moreover, if the basis $\mathbf{e}$ is in
addition shrinking, then its topological dual can be isometrically
identified with the space $(\bigoplus_i X_i^*)_{\mathbf{e^*}}$,
where $\mathbf{e^*}$ stands for the 1-unconditional basis of $E^*$
formed by the biorthogonal functionals associated
with~$\mathbf{e}$. When $\mathbf{e}$ is the unit vector basis of
$\ell_p$ $(1\leq p\leq\infty)$ we write $(\bigoplus_i X_i)_p$
instead of $(\bigoplus_i X_i)_{\mathbf{e}}$.

Given a Banach space, $E$, with a 1-unconditional basis
$\mathbf{e}=(e_i)$ we denote by $E^m$ the space $[e_i]_{i=1}^m$.
If $X$ is a Banach space we denote by $E^m(X)$ (resp.~$E(X)$) the
Banach space $(\bigoplus_i X_i)_{\mathbf{e}}$, where $X_i = X$,
$1\leq i\leq m$, and $X_i = \{0\}$, $i>m$ (resp.~$X_i = X$ for all
$i$). In particular, $\ell_p(X)$ (resp.~$\ell_p^m(X)$) denotes the
$\ell_p$-sum of countably infinitely many (resp.~$m$) copies
of~$X$. When appropriate, we may for $n\in \Bbb N$ identify
$E^n(X)$ with $E^n\otimes X$.

Given Banach spaces $X$ and $Y$ we write $\cal A(X,Y)$
(resp.~$\cal F(X,Y)$) for the Banach (resp.~normed) space of
approximable (resp.~finite-rank) operators from $X$ to $Y$. When
appropriate we shall identify $\cal F(X,Y)$ with $X^*\otimes Y$,
so that for $x^*\in X^*$, $y\in Y$ the rank-1 operator $x\mapsto
x^*(x)y$ is denoted $x^*\otimes y$. When $X = Y$ we simply write
$\cal A(X)$ (resp.~$\cal F(X)$). Likewise we shall use tensor
notation for operators $E^m(X)\mapsto E^n(X)$ for a Banach space
$E$ with a 1-unconditional basis and an arbitrary Banach space $X$
so that, for $m,n\in\Bbb N$ we identify $\cal A(E^m(X),E^n(X))$
with $\cal A(E^m,E^n)\otimes\cal A(X)$.

For any Banach space $X$ and positive integers $n>m$, there is a
natural isometric embedding $E^m(X)\hookrightarrow E^n(X)$ which
in turn induces an isometric Banach algebra homomorphism $\cal
A(E^m(X))\hookrightarrow \cal A(E^n(X))$. Letting $m$ and $n$ vary
we obtain a direct system of Banach algebras and isometric Banach
algebra homomorphisms. Its inductive limit is also a Banach
algebra that we denote by $\cal A_0(E(X))$. Note that $\cal
A_0(\ell_p(X)) = \cal A(\ell_p(X))$, $1<p\leq\infty$.

Recall that a Banach space $X$ is said to have the
$\lambda$-bounded approximation property ($\lambda$-BAP in short)
if there is a net $(T_\alpha)\subset\cal F(X)$ of bound $\lambda$
converging strongly to the identity operator on $X$. We write this
last as $T_\alpha\stackrel{s}{\to}I_X$. If, in addition, the
$T_\alpha$'s can be chosen to be projections then $X$ is called a
$\pi_\lambda$-space. A Banach space is said to have the bounded
approximation property (BAP in short) if it has the $\lambda$-BAP
for some $\lambda$, and is said to be a $\pi$-space if it is a
$\pi_\lambda$-space for some $\lambda$.

Recall that a bounded net $(e_\alpha)$ in a normed algebra $\cal
A$ is called a bounded approximate identity (BAI in short) for
$\cal A$ if $\lim_\alpha e_\alpha a = \lim_\alpha a e_\alpha =
a\,$ $(a\in\cal A)$. A normed $\cal A$-bimodule, $\cal X$ is {\it
essential}, if $\cal A\cdot\cal X\cdot\cal A$ is dense in $\cal
X$. Clearly, if $\cal A$ has a BAI, then this BAI is also a BAI
for any essential $\cal A$-bimodule. It is well known that the
algebra of approximable operators on a Banach space $X$ has a BAI
of bound $\lambda$ if and only if $X^*$ has the $\lambda$-BAP
\cite[Theorem~3.3]{GW}, \cite{Sa}.

Lastly, there is an intrinsic characterization of amenability that
is particularly useful in this setting. Precisely, a Banach
algebra $\cal A$ is amenable if and only if it has an approximate
diagonal, i.e., a bounded net $(d_\alpha)$ in $\cal
A\widehat{\otimes}\cal A$ such that $\pi(d_\alpha)a\to a$ and $a
d_\alpha - d_\alpha a\to 0$ $(a\in\cal A)$, where $\pi:\cal
A\widehat{\otimes}\cal A\to\cal A$, $a\otimes b\mapsto ab$
\cite[Lemma~1.2 and Theorem~1.3]{J2}. The Banach algebra $\cal A$
is said to be $K$-amenable if it has an approximate diagonal of
bound $K$. The smallest such $K$ is called the amenability
constant of $\cal A$.

Other definitions and results shall be given as they are needed.

\section{Property $\mathbb{A}$ revised}

Recall from \cite{GJW} that a Banach space $X$ is said to have
property $\Bbb A$ if there exist a constant $K>0$ and a bounded
net of projections $(P_\alpha)\subset\cal A(X)$ such that
\begin{enumerate}
\item[i)] $P_\alpha\stackrel{s}{\to} I_X$;

\item[ii)] $P_\alpha^*\stackrel{s}{\to} I_{X^*}$;

\item[iii)] For each $\alpha$ there is a finite group
$G_\alpha\subset \cal F(X_\alpha)$ whose linear span is $\cal
F(X_\alpha)$ and such that $\max_{T\in G_\alpha} \|T\| \leq K$
(where $X_\alpha = \text{rg}\,P_\alpha$).
\end{enumerate}

Property $\Bbb A$ was introduced in \cite{GJW} in an attempt to
explain amenability of $\cal A(X)$ as a consequence of some sort
of approximation property. Indeed, Banach spaces with this
property must carry amenable algebras of approximable operators
\cite[Theorem~4.2]{GJW}. Though we believe the converse is
unlikely to be true, no example of a Banach space $X$ without
property $\Bbb A$ and so that $\cal A(X)$ is amenable seems to be
known. The main result of this section, Corollary~\ref{gd} below,
is a characterization of amenability of the algebra of
approximable operators in terms of a property analogous to
property $\mathbb{A}$.

We start with the following.

\begin{prop}\label{am}
Let $X$ be a Banach space such that $\cal A(X)$ is $K$-amenable.
Suppose in addition that $\cal A(X)$ contains a bounded net of
projections, $(P_\alpha)_{\alpha\in A}$, such that $P_\alpha
\stackrel{s}{\to} I_X$ and $P_\alpha^*\stackrel{s}{\to} I_{X^*}$.
Then $\cal A(X)$ has an approximate diagonal
$(\delta_\alpha)_{\alpha\in A}$ with the following properties:
\begin{enumerate}
\item[$a)$] $\limsup_\alpha \|\delta_\alpha\|_{\wedge} \leq
\lambda K$, where $\lambda = \limsup_\alpha \|P_\alpha\|$;
\vspace{1mm}

\item[$b)$] $\pi(\delta_\alpha) = P_\alpha$ $(\alpha\in A)$;
\vspace{1mm}

\item[$c)$] $W\cdot\delta_\alpha = \delta_\alpha\cdot W$ for every
$W\in P_\alpha \cal A(X) P_\alpha$ $(\alpha\in A)$; and
\vspace{1mm}

\item[$d)$] For every $\alpha\in A$ there exists $\beta =
\beta(\alpha)\in A$ such that $\delta_\alpha\in \cal A(X) P_\beta
\otimes P_\beta \cal A(X)$.
\end{enumerate}
\end{prop}

\begin{proof}
Let $(d_i)_{i\in I}$ be an approximate diagonal for $\cal A(X)$
bounded by $K$. Since $P_\alpha \stackrel{s}{\to} I_X$ and
$P_\alpha^*\stackrel{s}{\to} I_{X^*}$, we can assume, without loss
of generality, that for every $i\in I$ there exists $\beta_i\in A$
such that $d_i\in\cal A(X) P_{\beta_i}\otimes P_{\beta_i} \cal
A(X)$. Let $x\in X$ and $x^*\in X^*$ be fixed vectors such that
$x^*(x) = 1$, and let $(\varepsilon_\alpha)$ be a net of positive
numbers converging to zero ($\varepsilon_\alpha = 1/\text{rk}\,
P_\alpha$ $(\alpha\in A)$ will do).

Given $i\in I$, let $\Phi_i: X^*\otimes X \to \cal F(X)\otimes\cal
F(X)$ be the linear map which is defined on elementary tensors by
$\Phi_i(\xi^*\otimes \xi) := x^*\otimes \xi\cdot
d_i\cdot\xi^*\otimes x$. It is readily seen that $\Phi_i$ is an
$\cal F(X)$-bimodule morphism.

For each $\alpha\in A$ choose $i(\alpha)\in I$ `big enough' so
that
\begin{equation}\label{1}
\left|1 - x^*\left(\pi\left(d_{i(\alpha)}\right)x\right)\right|
\leq \varepsilon_\alpha,
\end{equation}
and
\begin{equation}\label{2}
\left\|\Phi_{i(\alpha)}(P_\alpha) - P_\alpha\cdot
d_{i(\alpha)}\right\|_{\wedge} \leq \varepsilon_\alpha.
\end{equation}
Since $(\pi(d_i))$ is a bounded approximate identity for $\cal
A(X)$, it is clear that (\ref{1}) holds for every $i\in I$ `big
enough'. To see that the same is true about (\ref{2}) note that
for every $i\in I$ we have
\begin{eqnarray*}
\Phi_i(\xi^*\otimes \xi) - \xi^*\otimes\xi\cdot d_i &=&
x^*\otimes\xi\cdot d_i\cdot\xi^*\otimes x -
\xi^*\otimes\xi\cdot d_i \\
&=& x^*\otimes\xi\cdot \big(d_i\cdot\xi^*\otimes x - \xi^*\otimes
x\cdot d_i\big) \qquad (\xi\in X,\,\xi^*\in X^*).
\end{eqnarray*}
This last, combined in the obvious way with the facts that
$P_\alpha$ is finite-rank and that $(d_i)$ is an approximate
diagonal, gives the desired conclusion.

Now define a new net $(\delta_\alpha)\in\cal
A(X)\widehat{\otimes}\cal A(X)$ by
\[
\delta_\alpha := \varkappa_\alpha \Phi_{i(\alpha)}(P_\alpha)\qquad
(\alpha\in A),
\]
where $\varkappa_\alpha = 1/x^*(\pi(d_{i(\alpha)})x)$. We show
next that $(\delta_\alpha)$ has all required properties. First
note that
\begin{eqnarray*}
\|\delta_\alpha\|_{\wedge} &=& \varkappa_\alpha
\|\Phi_{i(\alpha)}(P_\alpha)\|_{\wedge} \;\leq\; \varkappa_\alpha
\|P_\alpha\cdot d_{i(\alpha)}\|_{\wedge} + \varkappa_\alpha
\|\Phi_{i(\alpha)}(P_\alpha) -
P_\alpha\cdot d_{i(\alpha)}\|_{\wedge} \\
&\leq& \varkappa_\alpha \|P_{\alpha}\| K + \varkappa_\alpha
\varepsilon_\alpha ,
\end{eqnarray*}
and so, $\limsup_\alpha \|\delta_\alpha\|_{\wedge} \leq \lambda
K$, that is, $(a)$ is satisfied.

That $(\delta_\alpha)$ satisfies $(b)$ follows immediately from
its definition above and the definition of $\Phi_i$. As for $(c)$,
just recall that $\Phi_{i(\alpha)}$ is an $\cal F(X)$-bimodule
morphism so $W\cdot\delta_\alpha = \delta_\alpha\cdot W$ whenever
$W P_\alpha = P_\alpha W$. By our assumption about $(d_i)$, at the
beginning of the proof, it is clear that $(d)$ is satisfied too.

Finally, since $P_\alpha\stackrel{s}{\to} I_X$ and $P_\alpha^*
\stackrel{s}{\to} I_{X^*}$, we have that
\[
\quad W\cdot\delta_\alpha - \delta_\alpha\cdot W = (W - P_\alpha W
P_\alpha)\cdot\delta_\alpha + \delta_\alpha\cdot(P_\alpha W
P_\alpha - W)\to 0 \qquad (W\in\cal A(X)).
\]
Obviously, $(\pi(\delta_\alpha))$ is a BAI for $\cal A(X)$, so,
$(\delta_\alpha)$ is an approximate diagonal for $\cal A(X)$.
\end{proof}

Thus, if $\cal A(X)$ is amenable and has a net of projections as
in the lemma, then it has an approximate diagonal whose elements
behave themselves like diagonals in a sense that we make more
precise in the next definition.

\begin{definition}
Let $X$ and $Y$ be finite-dimensional Banach spaces, and let $\cal
A$ be a subalgebra of $\cal F(X)$. We call an element $\Delta\in
\cal F(Y,X)\widehat{\otimes}\cal F(X,Y)$ a generalized diagonal
$($g.d. in short$)$ for $\cal A$, if
\begin{enumerate}
\item[i)] $W \Delta = \Delta W$ $(W\in\cal A)$; and \vspace{1mm}

\item[ii)] $\pi(\Delta)W = W$ $(W\in\cal A)$.
\end{enumerate}
\end{definition}

It is easily seen that when $\cal A = \cal F(X)$, an element
$\Delta\in\cal F(Y,X)\widehat{\otimes}\cal F(X,Y)$ is a
generalized diagonal for $\cal A$ if and only if there exists an
$\cal A$-bimodule morphism $\rho:\cal A\to\cal
F(Y,X)\widehat{\otimes}\cal F(X,Y)$ so that $\pi \circ \rho =
I_{\cal A}$ and $\rho(I_X) = \Delta$. Furthermore, if
$(x_k)_{k=1}^m$ and $(y_i)_{i=1}^n$ are bases of $X$ and $Y$,
respectively, then it follows from this last observation, that
$\Delta$ can be written as
\begin{equation}\label{gf}
\Delta = \sum\nolimits_{i,j} a_{i,j}\,\sum\nolimits_k
(y_j^*\otimes x_k)\otimes (x_k^*\otimes y_i),
\end{equation}
for some scalars $a_{i,j}$ satisfying $\sum_i a_{i,i} = 1$, where,
as is customary, the $y_j^*$'s (resp.~the $x_k^*$'s) denote the
biorthogonal functionals associated with the basis $(y_i)_{i=1}^n$
(resp.~$(x_k)_{k=1}^m$). Conversely, it can be easily verified
that every element of the form (\ref{gf}) is a generalized
diagonal for~$\cal F(X)$.

Now the main result of this section is merely a restatement of
Proposition~\ref{am} in terms of generalized diagonals.

\begin{cor}\label{gd}
Let $X$ be a Banach space. Suppose $\cal A(X)$ contains a bounded
net of projections, $(P_\alpha)_{\alpha\in A}$, such that
$P_\alpha\stackrel{s}{\to} I_X$ and $P_\alpha^*\stackrel{s}{\to}
I_{X^*}$. Set $X_\alpha = \mathrm{rg}\,P_\alpha$ $(\alpha\in A)$.
Then $\cal A(X)$ is amenable if and only if there is a constant
$K>0$ such that for every $\alpha\in A$ there exists $\beta =
\beta(\alpha)\in A$ such that $\cal F(X_\beta,X_\alpha)
\widehat{\otimes}\cal F(X_\alpha,X_\beta)$ contains a generalized
diagonal for $\cal F(X_\alpha)$ of norm no greater than $K$.
\end{cor}

\begin{proof}
First suppose $\cal A(X)$ is amenable. By Proposition~\ref{am},
$\cal A(X)$ has an approximate diagonal,
$(\delta_\alpha)_{\alpha\in A}$, satisfying $(a)$--$(d)$ of the
same proposition. For each $\alpha\in A$, let $\beta =
\beta(\alpha)\in A$ be as in $(d)$. Let $P_\alpha^c$ (resp.
$P_\beta^c$) denote the corestriction of $P_\alpha$ (resp.
$P_\beta$) to its range, and let $\imath_\alpha:X_\alpha\to X$
(resp. $\imath_\beta:X_\beta\to X$) denote the canonical embedding
of $X_\alpha$ (resp. $X_\beta$) into $X$. Then define
$\Delta_\alpha\in\cal F(X_\beta,X_\alpha) \widehat{\otimes}\cal
F(X_\alpha,X_\beta)$ by $\Delta_\alpha := \Phi_\alpha
(\delta_\alpha)$, where $\Phi_{\alpha}:\cal
A(X)\widehat{\otimes}\cal A(X)\to\cal
F(X_\beta,X_\alpha)\widehat{\otimes}\cal F(X_\alpha,X_\beta)$ is
the linear map defined on elementary tensors by $\Phi_\alpha
(R\otimes S) := P_\alpha^c R\, \imath_\beta \otimes P_\beta^c S\,
\imath_\alpha$ $(R,S\in\cal A(X))$. It is easy to verify that
$\Delta_\alpha$ is a g.d. for $\cal F(X_\alpha)$ $(\alpha\in A)$.
The desired conclusion now follows on noting that the family
$(\Phi_\alpha)_{\alpha\in A}$ is uniformly~bounded.

Conversely, for each $\alpha\in A$, let $\Delta_\alpha\in\cal
F(X_\beta,X_\alpha)\widehat{\otimes}\cal F(X_\alpha,X_\beta)$ be a
g.d. for $\cal F(X_\alpha)$ of norm $\,\leq K$, and let
$\Psi_{\alpha}:\cal F(X_\beta,X_\alpha)\widehat{\otimes}\cal
F(X_\alpha,X_\beta)\to\cal A(X)\widehat{\otimes}\cal A(X)$ be the
linear map defined on elementary tensors by $\Psi_\alpha (U\otimes
V) := \imath_\alpha U P_\beta^c \otimes \imath_\beta V P_\alpha^c$
$(U\in\cal F(X_\beta,X_\alpha),\,V\in\cal F(X_\alpha,X_\beta))$.
Then $(\Psi_\alpha(\Delta_\alpha))$ is an approximate diagonal for
$\cal A(X)$.
\end{proof}

\begin{rem}
If $X$ is not a $\pi$-space but still $X^*$ has the BAP, as must
be the case if $\cal A(X)$ is amenable \cite{GJW}, then we can
argue as follows. First, we choose a net of projections
$(P_\alpha)$ in $\cal F(X)$ such that $P_\alpha\stackrel{s}{\to}
I_X$ and $P_\alpha^* \stackrel{s}{\to} I_{X^*}$. Such a net, of
course, would be necessarily unbounded. Then we choose a bounded
net $(T_\alpha)$ in $\cal F(X)$ such that $\,P_\alpha T_\alpha =
P_\alpha = T_\alpha P_\alpha\,$ for every $\alpha$, and set
$X_\alpha := \text{rg}\,T_\alpha$. It can be shown that $\cal
A(X)$ is amenable if and only if there is a constant $K>0$ such
that for every $\alpha\in A$ there exists $\beta =
\beta(\alpha)\in A$ such that $\,\cal
F(X_\beta,X_\alpha)\widehat{\otimes}\cal F(X_\alpha,X_\beta)\,$
contains a generalized diagonal for $\cal A_\alpha =
P_\alpha|^{X_\alpha}\cal A(X) P_\alpha|_{X_\alpha}(\subseteq\cal
F(X_\alpha))$ of norm no greater than $K$ (here
$P_\alpha|_{X_\alpha}$ and $P_\alpha|^{X_\alpha}$ denote the
restriction and corestriction, respectively, of $P_\alpha$ to
$X_\alpha$).
\end{rem}

\begin{rem}
Note that (iii) of the definition of property $\Bbb A$ guarantees
the existence of a diagonal (and hence a generalized diagonal) for
$\cal F(X_\alpha)$ in $\cal F(X_\alpha)\widehat{\otimes}\cal
F(X_\alpha)$ whose norm does not exceed $K$, namely,
$\frac{1}{|G_\alpha|} \sum_{T\in G_\alpha} T\otimes T^{-1}$.
\end{rem}

\begin{ex}\label{lpq}
Let $(n_k)$ be an unbounded sequence of positive integers, and let
$1\leq p\ne q\leq\infty$. It is shown in \cite[Theorem~6.5]{GJW}
that the algebra $\cal A\big((\bigoplus_k \ell_p^{n_k})_q\big)$ is
amenable. It seems to be unknown whether or not this algebra has
property $\mathbb{A}$. However, it is relatively easy to show that
this algebra satisfies the condition of Corollary~\ref{gd}.
Indeed, fix $i\in\Bbb{N}$ and let $m = \max\{n_1,\ldots,n_i,i\}$.
The algebra $\cal F\big(\ell_q^m(\ell_p^m)\big)$ has a diagonal
$\Delta_m$ of norm 1 (see the discussion below). (Furthermore,
note that $\Delta_m$ can be given explicitly.) As $(n_k)$ is
unbounded, there are positive integers $k_1<k_2<\ldots<k_m$ so
that $m\leq\min\{n_{k_j}:1\leq j\leq m\}$. Clearly, we can think
of $(\bigoplus_{j=1}^i \ell_p^{n_j})_q =: X_i$
(resp.~$\ell_q^m(\ell_p^m)$) as a 1-complemented subspace of
$\ell_q^m(\ell_p^m)$ (resp.~$(\bigoplus_{j=1}^{k_m}
\ell_p^{n_j})_q =: X_{k_m}$). Let $P_1:\ell_q^m(\ell_p^m)\to X_i$
and $P_2:X_{k_m}\to\ell_q^m(\ell_p^m)$ be the natural projections,
and let $\imath_1:X_i\to\ell_q^m(\ell_p^m)$ and
$\imath_2:\ell_q^m(\ell_p^m)\to X_{k_m}$ be the corresponding
inclusion maps. It is easy to see that the image of $\Delta_m$ by
the linear map $\cal F(\ell_q^m(\ell_p^m))\widehat{\otimes}\cal
F(\ell_q^m(\ell_p^m))\to \cal F(X_{k_m},X_i)\widehat{\otimes}\cal
F(X_i,X_{k_m})$, $R\otimes S\mapsto
P_1RP_2\otimes\imath_2S\imath_1$, is a generalized diagonal for
$\cal F(X_i)$ in $\cal F(X_{k_m},X_i)\widehat{\otimes}\cal
F(X_i,X_{k_m})$ of norm at most 1. 
The rest is clear.
\end{ex}

It can be shown that if $X$ is a Banach space so that $\cal A(X)$
is $K$-amenable then $\cal A(\ell_p^n(X))$ is $K$-amenable for
every $1\leq p\leq\infty$ and $n\in\Bbb N$. Indeed, let $\cal H$
be the group of permutation matrices generated by a cyclic
permutation of the unit vector basis of $\ell_p^n$, and let $\cal
G = \{\text{diag}(t)\sigma : t\in\{\pm 1\}^n, \sigma\in\cal H\}$,
so $\frac{1}{|\cal G|}\sum_{g\in\cal G} g\otimes g^{-1}$ is a
diagonal for $\cal A(\ell_p^n)$ \cite[Example~3.3]{GJW}. Let
$(d_\alpha)$ be an approximate diagonal for $\cal A(X)$ of bound
$K$, and choose for each $d_\alpha$ a representation $\sum_j
U_{\alpha,j}\otimes V_{\alpha,j}$ such that $\sum_j
\|U_{\alpha,j}\|\|V_{\alpha,j}\| \leq K$. Then the elements
$\delta_\alpha := \frac{1}{|\cal G|}\sum_{j,g} (g\otimes
U_{\alpha,j})\otimes(g^{-1}\otimes V_{\alpha,j})\in\cal
A(\ell_p^n(X))\widehat{\otimes}\cal A(\ell_p^n(X))$ form an
approximate diagonal for $\cal A(\ell_p^n(X))$ of bound $K$.
Crucial in establishing this last is the fact that the $g$'s are
permutation matrices, as it seems $\ell_p(X)$ is rarely ever a
tight tensor product in the sense of \cite[Definition~2.1]{GJW}.
This is better exemplified through our next result, which extends
Theorem~2.5 of \cite{GJW}.

\begin{prop}\label{tight}
Let $E$ be a Banach space with a 1-unconditional basis $\mathbf{e}
= (e_n)$. $($Recall $E^n = [e_i]_1^n$.$)$ Suppose there is $K>0$
so that for each $m\in\Bbb N$ there exists $n\geq m$ such that
$\cal F(E^n,E^m)\widehat{\otimes}\cal F(E^m,E^n)$ contains a
generalized diagonal $\Delta_m$ with the following property: there
exists a representation $\sum_{i=1}^k R_{m,i}\otimes S_{m,i}$ of
$\Delta_m$ such that $\sum_i \|R_{m,i}\|\|S_{m,i}\|\leq K$ and the
matrix representation of each $R_{m,i}$ $($resp.~$S_{m,i})$ with
respect to the $e_i$'s has at most one non-zero entry in each row
and column. If $X$ is a Banach space such that $\cal A(X)$ is
$M$-amenable then $\cal A_0(E(X))$ is $KM$-amenable.
\end{prop}

\begin{proof}
Let $(d_\alpha)$ be an approximate diagonal for $\cal A(X)$ of
bound $M$. For each $d_\alpha$ choose a representation $\sum_j
U_{\alpha,j}\otimes V_{\alpha,j}$ such that $\sum_j
\|U_{\alpha,j}\|\|V_{\alpha,j}\|\leq M$. We show that the elements
$\delta_{m,\alpha} := \sum_{i,j} (R_{m,i}\otimes
U_{\alpha,j})\otimes(S_{m,i}\otimes V_{\alpha,j})$ form an
approximate diagonal of bound $KM$ for $\cal A_0(E(X))$. First
note that $R_{m,i}\otimes U_{\alpha,j}\in\cal
A(E^n(X),E^m(X))\subset\cal A_0(E(X))$ and for any set of vectors
$x_1,x_2,\ldots,x_n$ in $X$ we have
\begin{eqnarray*}
\left\|R_{m,i}\otimes U_{\alpha,j}\Big(\sum\nolimits_{k=1}^n
e_k\otimes x_k\Big)\right\| &=& \left\|\sum\nolimits_k
R_{m,i}(e_k)\otimes U_{\alpha,j}(x_k)\right\| \\
&=& \left\|\sum\nolimits_k \big\|U_{\alpha,j}(x_k)\big\|
R_{m,i}(e_k)\right\| \\
&\leq& \|R_{m,i}\|\|U_{\alpha,j}\|\left\|\sum\nolimits_k
e_k\otimes x_k\right\|,
\end{eqnarray*}
where the second equality follows from the fact that there is at
most one non-zero entry in each row and column of the matrix
representation of $R_{m,i}$ with respect to the $e_i$'s. Thus
$\|R_{m,i}\otimes U_{\alpha,j}\|\leq\|R_{m,i}\|\|U_{\alpha,j}\|$.
Likewise, $S_{m,i}\otimes V_{\alpha,j}\in\cal
A(E^m(X),E^n(X))\subset\cal A_0(E(X))$ and $\|S_{m,i}\otimes
V_{\alpha,j}\|\leq\|S_{m,i}\|\|V_{\alpha,j}\|$. Combining these
estimates we readily obtain that
$\|\delta_{m,\alpha}\|_{\wedge}\leq KM$.

In order to verify that $\pi(\delta_{m,\alpha})W\to W$ and
$W\cdot\delta_{m,\alpha} - \delta_{m,\alpha}\cdot W\to 0$
$\,(W\in\cal A_0(E(X))$ it is clearly enough to look at operators
$W$ of the form $E_{rs}\otimes T$ where $E_{rs} = e_r^*\otimes
e_s$ and $T\in\cal A(X)$. The procedures are standard, so we leave
the details to the reader.
\end{proof}

An immediate consequence of the above is the following.

\begin{cor}
Let $(n_k)$ be an increasing sequence of positive integers, let
$1\leq p,q\leq\infty$, and let $X$ be a Banach space such that
$\cal A(X)$ is amenable. Then $\cal A_0\big((\bigoplus_k
\ell_p^{n_k}(X))_q\big)$ is amenable. In particular, if $q>1$ then
$\cal A\big((\bigoplus_k \ell_p^{n_k}(X))_q\big)$ is amenable.
\end{cor}

\begin{proof}
The space $(\bigoplus_k \ell_p^{n_k})_q$ satisfies all hypotheses
of Theorem~\ref{tight} (see Example~\ref{lpq} above).
\end{proof}

Corollary~\ref{gd} essentially reduces the study of amenability of
algebras of approximable operators on $\pi$-spaces to the problem
of finding the minimum among the norms of all generalized
diagonals for $\cal F(X)$ in $\cal F(Y,X)\widehat{\otimes}\cal
F(X,Y)$ with $X$ and $Y$ finite-dimensional. Here, of course, the
main difficulty arises in estimating the projective norm. In some
cases, this task can be further simplified. For instance, let the
basis $(y_i)_{i=1}^n$ of $Y$ be 1-unconditional. Set $p_{i,j} =
\sum_k (y_j^*\otimes x_k)\otimes(x_k^*\otimes y_i)$ $\,(1\leq
i,j\leq n)$. Then, while looking for generalized diagonals of
minimum norm, we can restrict our attention to convex linear
combinations of the $p_{i,i}$'s. Indeed, in this case we have that
\begin{equation}\label{diagu}
\Big\|\sum\nolimits_{i,j} a_{i,j}\, p_{j,i}\Big\|_{\wedge} \geq
\Big\|\sum\nolimits_i a_{i,i}\, p_{i,i}\Big\|_{\wedge} =
\Big\|\sum\nolimits_i |a_{i,i}|\, p_{i,i}\Big\|_{\wedge}.
\end{equation}
To see this, consider the linear operator
\[
\Phi:\cal F(Y,X)\Hat{\otimes}\cal F(X,Y) \to \cal
F(Y,X)\Hat{\otimes}\cal F(X,Y),\; R\otimes S \mapsto 2^{-n}
\sum\nolimits_{t\in\{-1,1\}^n} R U_t\otimes U_t S,
\]
where $U_t\in\cal F(Y)$ is defined by $U_t(y_j) := t_j y_j$
$(1\leq j\leq n)$. It is clear that $\|\Phi\|\leq 1$, and it is
not difficult to see that $\Phi(\sum_{i,j} a_{i,j}\, p_{j,i}) =
\sum_i a_{i,i}\, p_{i,i}$, whence the inequality. As for the
equality, let $\Lambda\in\cal F(Y)$ be defined by $\Lambda(y_i) :=
\lambda_i y_i$, where $\lambda_i = \overline{a_{i,i}}/|a_{i,i}|$
$(1\leq i\leq n)$, and let $\Phi_\Lambda$ be the linear map
defined by
\[
R\otimes S \mapsto R\otimes \Lambda S \qquad (R\in\cal F(Y,X),\,
S\in\cal F(X,Y)).
\]
Then $\Phi_\Lambda$ is an isometry and
$\Phi_\Lambda(a_{i,i}\,p_{i,i}) = |a_{i,i}|\,p_{i,i}$ $(1\leq
i\leq n)$, so the equality follows. The claim that we can restrict
our attention to `convex' linear combinations now follows on
combining (\ref{diagu}) with the fact that the sum of the diagonal
coefficients in the representation $(\ref{gf})$ must be~1.

\begin{rem}
It is not hard to see that the sequence $(p_{i,i})$ has the same
basis, unconditional and symmetric constants as the basis $(y_i)$.
\end{rem}

It was asked in \cite{GJW} whether or not the $C_p$ spaces of
W.~B.~Johnson ($1<p<\infty$) carry amenable algebras of compact
operators. This question has an interesting interpretation in
terms of generalized diagonals. We consider the following more
general situation.

Let $(X_n)$ be a sequence of finite-dimensional Banach spaces
dense in the Banach--Mazur sense in the class of all
finite-dimensional Banach spaces, and let $\mathbf{e}$ be an
unconditional shrinking Schauder basis. Define $C_{\mathbf{e}} :=
(\bigoplus_n X_n)_\mathbf{e}$. It is readily seen from
Corollary~\ref{gd} that the algebra $\cal A(C_{\mathbf{e}})$ is
amenable if and only if there exists an absolute constant $K$ with
the following property: \vspace{2mm}

\noindent{\it For every finite-dimensional Banach space $X$ there
exists a finite-dimensional Banach space $Y$ so that $\cal
F(Y,X)\widehat{\otimes}\cal F(X,Y)$ contains a generalized
diagonal for $\cal F(X)$ of norm at most $K$.} \vspace{2mm}

We do not know if one such constant can exist. However, if $X$ is
a finite-dimensional Banach space with unconditional constant
$<\lambda$ then, by a finite-dimensional version of a well known
result of J.~Lindenstrauss \cite[Remark~4]{L}, there exists a
finite-dimensional Banach space $Y$ with symmetric constant
$<\lambda$ such that $X$ is a 1-complemented subspace of $Y$.
Thus, $\cal F(Y)\widehat{\otimes}\cal F(Y)$ contains a diagonal
for $\cal F(Y)$ of norm $<\lambda$, and in turn $\cal
F(Y,X)\widehat{\otimes}\cal F(X,Y)$ contains a generalized
diagonal for $\cal F(X)$ of norm $<\lambda$. As a simple
consequence of this we quote the following.

\begin{prop}
Let $(X_n)$ be a sequence of finite-dimensional Banach spaces with
unconditional constant $<\lambda$, dense in the Banach-Mazur sense
in the class of all finite-dimensional Banach spaces with
unconditional constant $<\lambda$. Let $\mathbf{e}$ be an
shrinking 1-unconditional Schauder basis. Then the algebra $\cal
A\big((\bigoplus_n X_n)_\mathbf{e}\big)$ has property
$\mathbb{A}$.
\end{prop}

We should like to end this section by noting that, there is an
analogue of Lindenstrauss's result, due to Johnson, Rosenthal and
Zippin, which states that there is a universal constant $C~~(\leq
16^{12})$ so that every finite-dimensional Banach space is a
1-complemented subspace of a finite-dimensional space with basis
constant no greater than $C$ \cite[Corollary~4.12(a)]{JRZ}.

\section{Amenability and \\ equivalence of operator ideal norms}

Unfortunately, the characterization found in the previous section
does not help very much when it comes to determine if the algebra
of approximable operators on a given Banach space is amenable or
not. In this section we take a different approach.

Recall that a Banach space $X$ is called \emph{approximately
primary} if, for every projection $P\in\cal B(X)$, at least one of
the product maps $\pi:\cal A(PX,X)\widehat{\otimes}\cal
A(X,PX)\to\cal A(X)$ or $\pi:\cal A((I-P)X,X)\widehat{\otimes}\cal
A(X,(I-P)X)\to\cal A(X)\,$ is surjective. This notion was
introduced in \cite{GJW}, where it was shown that if $\cal A(X)$
is amenable then $X$ must be approximately primary. Moreover, also
in the same paper (see the proof of \cite[Theorem~6.9]{GJW} and
comments after it), it was shown that none of the following spaces
is approximately primary: $\ell_p\oplus \ell_q$ $(1<p,q<\infty$,
$p\ne q$ and neither $p$ nor $q$ is equal 2), $\ell_1\oplus
\ell_q$ $(q>2)$ and $\ell_p\oplus \ell_\infty$ $(p<2)$.

Essential to the proof of this last result were the following:

\newcounter{q}
\begin{list}{\arabic{q}.}{\usecounter{q}
\setlength{\leftmargin}{2cm}}

\item[Fact~1)] $\;$Given Banach spaces $X$ and $Y$, if the space
$X\oplus Y$ is approximately primary, then at least one of the
product maps $\pi:\cal A(Y,X)\widehat{\otimes}\cal A(X,Y)\to\cal
A(X)$ or $\pi:\cal A(X,Y)\widehat{\otimes}\cal A(Y,X)\to\cal
A(Y)\,$ is surjective.

\item[Fact~2)] $\;$For every Banach space $X$, the product map
$\pi:\cal A(\ell_p,X)\widehat{\otimes}\cal A(X,\ell_p)\to\cal
A(X)$ is surjective if and only if the bilinear map $\varphi:\cal
A(\ell_p,X)\times\cal A(X,\ell_p)\to\cal A(X)$ is open.
\end{list}

The results of this section are, to some extent, generalizations
of these two facts. We start by recalling some standard
terminology.

Let $\cal F$ be the operator ideal of all finite-rank operators
between Banach spaces so, for every pair of Banach spaces $(X,Y)$
we have $\cal F\cap\cal B(X,Y) = \cal F(X,Y)$. Recall that an
\emph{operator ideal norm on $\cal F$} is a function $\gamma:\cal
F \to [0,\infty[$ that satisfies:
\begin{enumerate}
\item[$a)$] $\gamma\big|_{\cal F(E,F)}$ is a norm for every pair
of Banach spaces $E$ and $F$; \vspace{1mm}

\item[$b)$] $\gamma(I_{\mathbb{C}}:\mathbb{C}\to \mathbb{C}) = 1$;
\vspace{1mm}

\item[$c)$] If $A\in\cal B(Y,Y_0) $, $B\in\cal B(X_0,X)$ and
$T\in\cal F(X,Y)$ then $\gamma(ATB) \leq \|A\|\gamma(T)\|B\|$.
\end{enumerate}
It is well known that if $\gamma$ is as above then
$\|T\|\leq\gamma(T)$ for every $T\in\cal F$. Moreover, if $T =
f\otimes x$ then $\|T\| = \|f\|\|x\| = \gamma(T)$ $(x\in X,f\in
X^*)$.

In the terminology of \cite[\S 9]{DF} the operator ideal $\cal F$
endowed with an operator ideal norm as in the above definition is
a normed operator ideal. Of course, it will not be `Banach' (i.e.,
complete). The reason for doing things in this way should become
clear later on. Examples of operator ideal norms on $\cal F$ are
the restrictions of the classical operator ideal norms, like
nuclear and $\pi$-summing norms, to $\cal F$.

Now the main result of this section reads as follows.

\begin{thm}\label{2n}
Let $X$ and $Y$ be Banach spaces, and let $\gamma$ and $\tau$ be
operator ideal norms on $\cal F$. Suppose that
\begin{enumerate}
\item[i)] $\cal A(X)$ has a BAI of bound $\lambda$;

\item[ii)] The multiplication $\pi:\cal
A(Y,X)\widehat{\otimes}\cal A(X,Y)\to\cal A(X)$ is surjective with
inversion constant $\beta$;

\item[iii)] $\gamma$ and $\tau$ are equivalent on one of $\cal
F(Y,X)$ or $\cal F(X,Y)$, say $c\gamma\leq\tau\leq C\gamma$.
\end{enumerate}
Then $\gamma$ and $\tau$ are equivalent on $\cal F(X)$,
specifically
\[
c\beta^{-2}\lambda^{-2}\gamma\leq\tau\leq C\beta^2\lambda^2\gamma.
\]
\end{thm}

\begin{proof}
Let $F\in\cal F(X)$ and let $(T_\alpha)$ be a BAI for $\cal A(X)$ of
bound $\lambda$. Note that for any operator ideal norm $\gamma$ we
have that $\lim_\alpha \gamma(F-T_\alpha F T_\alpha)= 0$. Indeed,
simply write $F=GH$ with $G,H\in\cal F(X)$. Then
\begin{eqnarray*}
\gamma(F-T_\alpha FT_\alpha) &=&
\gamma((G - T_\alpha G)H + T_\alpha G(H - H T_\alpha)) \\
&\leq& \|G - T_\alpha G\| \gamma(H) + \lambda \gamma(G)\|H - H
T_\alpha\|,
\end{eqnarray*}
which tends to 0 as $\alpha\to\infty$.

Let $G\in\cal F(X)$ with $\|G\|\leq\lambda$ and let $L>\beta$.
Choose $\sum_i R_i\otimes S_i\in\cal A(Y,X)\widehat{\otimes}\cal
A(X,Y)$ so that $\sum_i R_i S_i = G$ and $\sum_i \|R_i\|\|S_i\| <
\lambda L$. Then for any $F\in\cal F(X)$ we have that
\[
\sum\nolimits_{i,j} \gamma(R_i S_i F R_j S_j) \leq \gamma(F)
\sum\nolimits_{i,j} \|R_i S_i\|\|R_j S_j\| \leq \gamma(F) L^2
\lambda^2.
\]
Moreover,
\begin{eqnarray*}
\lefteqn{\gamma\Big(\sum\nolimits_{1\leq i,j\leq n} R_i S_i F R_j
S_j - GFG\Big)} \hspace{2cm} \\ &=&
\gamma\Big(\Big(\sum\nolimits_1^n R_i S_i -
G\Big)F\Big(\sum\nolimits_1^n R_j S_j\Big) +
GF \Big(\sum\nolimits_1^n R_j S_j - G\Big)\Big) \\
&\leq& \Big\|\sum\nolimits_1^n R_i S_i - G\Big\|\gamma(F)\lambda L
+ \|G\|\gamma(F)\Big\|\sum\nolimits_1^n R_j S_j - G\Big\|,
\end{eqnarray*}
which tends to 0 as $n\to\infty$. So, the series $\sum_{i,j} R_i S_i
F R_j S_j$ is unconditionally $\gamma$-convergent in $\cal F(X)$ to
the sum $GFG$.

Assume that $\gamma$ and $\tau$ are equivalent on $\cal F(X,Y)$.
Then
\begin{eqnarray*}
c\gamma(GFG) &\leq& \sum\nolimits_{i,j} c\gamma(R_i S_i F R_j S_j)
\;\leq\; \sum\nolimits_{i,j} c\|R_i\|\gamma(S_i F)\|R_j S_j\| \\
&\leq& \sum\nolimits_{i,j} \|R_i\|\|S_i\|\tau(F)\|R_j S_j\|
\;\leq\; L^2 \lambda^2 \tau(F).
\end{eqnarray*}
Letting $G=T_\alpha$ and $\alpha\to\infty$ we obtain $c\gamma(F)
\leq L^2\lambda^2\tau(F)$. Likewise $\tau(F)\leq L^2 \lambda^2
C\gamma(F)$. A similar proof working with $F R_j$ rather than $S_i
F$ gives the result in case $\gamma$ and $\tau$ are equivalent on
$\cal F(Y,X)$.
\end{proof}

We now bring amenability into the picture. We start with the
following refinement of \cite[Theorem~6.8]{GJW}.

\begin{prop}\label{am+ap}
Let $X$ be a Banach space and let $P:X\to X$ be a bounded
projection. Set $Y=\mathrm{rg}\,P$ and $Z=\mathrm{rg}\,(I-P)$. If
$\cal A(X)$ is $K$-amenable then at least one of the maps
$\pi_1:\cal A(Z,Y)\widehat{\otimes}\cal A(Y,Z)\to\cal A(Y)$ or
$\pi_2:\cal A(Y,Z)\widehat{\otimes}\cal A(Z,Y)\to\cal A(Z)$ is
surjective with inversion constant no greater than
$4K\|P\|\|I-P\|\max\{\|P\|^3,\|I-P\|^3\}$.
\end{prop}

\begin{proof}
The proof is almost the same as that of \cite[Theorem~6.8]{GJW},
one only needs to keep track of the constants.

Set $P_1 = P$, $P_2 = I-P$, $\cal A = \cal A(X)$ and $\cal A_{ij}
= P_i\cal AP_j$ $\,(i,j=1,2)$. Let $A_{ii}^\circ = \pi(\cal
A_{ji}\widehat{\otimes}\cal A_{ij})$ with the norm
$\|\cdot\|^\circ$ inherited from $\cal A_{ji}\widehat{\otimes}\cal
A_{ij}/(\cal A_{ji}\widehat{\otimes}\cal A_{ij}\cap\ker{\pi})$ via
the natural isomorphism induced by the product map $\pi$
$\,(i,j=1,2,\,i\ne j)$. It is easy to see that $\|a_{ii}
a_{ii}^\circ\|^\circ\leq \|a_{ii}\|\|a_{ii}^\circ\|^\circ$
$(a_{ii}\in\cal A_{ii},\,a_{ii}^\circ\in\cal A_{ii}^\circ)$,
$i=1,2$. So $\cal A_{ii}^\circ$ is a Banach $\cal
A_{ii}$-bimodule, and $\cal A_{ii}$ is a Banach $\cal
A_{ii}^\circ$-bimodule, $i=1,2$. Let $\cal A^\circ = \{a\in\cal A
: P_iaP_i\in\cal A_{ii}^\circ,\, i=1,2\}$ with the norm
$\|a\|^\circ =
\max\{\|P_1aP_1\|^\circ,\|P_1aP_2\|,\|P_2aP_1\|,\|P_2aP_2\|^\circ\}$
$(a\in\cal A^\circ)$. Then $\|a a^\circ\|^\circ\leq
M\|a\|\|a^\circ\|^\circ$ and $\|a^\circ a\|^\circ\leq
M\|a\|\|a^\circ\|^\circ$ $(a\in\cal A,\,a^\circ\in\cal A^\circ)$
for some constant $M\leq 2\max\{\|P_1\|^2,\|P_2\|^2\}$, so $(\cal
A^\circ,\|\cdot\|^\circ)$ is a Banach $\cal A$-bimodule.

The map $D:\cal A\to\cal A^\circ$, $a\mapsto P_1aP_2 - P_2aP_1 =
P_1a - aP_1$ is a bounded derivation, and so, there is $C\in(\cal
A^\circ)^{**}$ such that $Da = aC - Ca$ $(a\in\cal A)$.
Furthermore, we can choose $C$ so that $\|C\|^\circ\leq KM\|D\|$
(see \cite[Theorem~1.3]{J2}). Let $C_{ii} = P_iCP_i\in(\cal
A_{ii}^\circ)^{**}$ and let $\imath_i:\cal A^\circ_{ii}\to\cal
A_{ii}$ be the inclusion map $(i=1,2)$. It can be shown that
$a_{ii}(\imath_i^{**}C_{ii}) = \lambda_i a_{ii} =
(\imath_i^{**}C_{ii}) a_{ii}$ $(a_{ii}\in\cal A_{ii})$ for some
$\lambda_i\in\Bbb C$ $(i=1,2)$. Moreover, $\lambda_2 - \lambda_1 =
1$ and if $\lambda_i\ne 0$ then $\cal A_{ii}^\circ = \cal A_{ii}$
and $\imath_i^{**}C_{ii} = C_{ii}$. (See the proof of
\cite[Theorem~6.8]{GJW} for details.)

As $\lambda_2 - \lambda_1 = 1$, at least one of $\lambda_1$ or
$\lambda_2$ must have absolute value greater than or equal $1/2$.
Without loss of generality suppose $|\lambda_1|\geq 1/2$, so we
have $\cal A_{11}^\circ = \cal A_{11}$ and $\imath_1^{**}C_{11} =
C_{11}$. Let $(e_\alpha)$ be a net in $\cal A_{11}^\circ$ bounded
by $\|\lambda_1^{-1} C_{11}\|^\circ$ and weak-$*$ convergent to
$\lambda_1^{-1} C_{11}$. Since $a_{11}^\circ(\lambda_1^{-1}C_{11})
= a_{11}^\circ = (\lambda_1^{-1}C_{11})a_{11}^\circ$, it is
readily seen that $e_\alpha a_{11}^\circ\to a_{11}^\circ$ and
$a_{11}^\circ e_\alpha\to a_{11}^\circ$ weakly for every
$a_{11}^\circ\in\cal A_{11}^\circ$. A standard argument (see
\cite[Proposition~2.9.14 (iii)]{Da}) shows that $\cal
A_{11}^\circ$ has BAI of bound $\|\lambda_1^{-1} C_{11}\|^\circ$.

Now let $a_{11}\in\cal A_{11}$ be arbitrary. It is easy to see
that $\cal A_{11}$ is an essential $\cal A_{11}^\circ$-bimodule,
so, by \cite[Theorem~2.9.24]{Da}, there exist $e^\circ\in\cal
A_{11}^\circ$ and $b\in\cal A_{11}$ such that $a_{11} = e^\circ
b$, $\|e^\circ\|^\circ\leq\|\lambda_1^{-1} C_{11}\|^\circ$ and
$\|b\|\leq \|a_{11}\|$. Thus,
\[
\|a_{11}\|^\circ\leq\|e^\circ\|^\circ\|b\|\leq
\|\lambda_1^{-1}C_{11}\|^\circ \|a_{11}\| \leq
2\|C\|^\circ\|a_{11}\|\leq 2KM\|D\|\|a_{11}\|.
\]
As $\|P_1aP_2 - P_2aP_1\|^\circ = \max\{\|P_1aP_2\|,\|P_2aP_1\|\}
\leq \|P_1\|\|P_2\|\|a\|$ $(a\in\cal A)$, we find that $\|D\| \leq
\|P_1\|\|P_2\|$, so
\[
\|a_{11}\|^\circ\leq 2KM\|P_1\|\|P_2\|\|a_{11}\|\qquad
(a_{11}\in\cal A_{11}).
\]

To finish the proof of the proposition, one just needs to note
that the linear isomorphisms $\cal A_{21}\widehat{\otimes}\cal
A_{12}\to\cal A(Z,Y)\widehat{\otimes}\cal A(Y,Z)$, $R\otimes
S\mapsto R|_Z^Y\otimes S|_Y^Z$, and $\cal A(Y)\to\cal A_{11}$,
$T\mapsto\imath T P_1$, where $\imath:Y\to X$ denotes the
inclusion map, have norms no greater than 1 and $\|P_1\|$,
respectively. Combining these two last estimates with those
previously found, we finally obtain that the inversion constant of
$\pi_1$ cannot be greater than $2KM\|P_1\|\|P_2\|
\max\{\|P_1\|,\|P_2\|\}$, as~claimed.
\end{proof}

Combining Proposition~\ref{am+ap} and Theorem~\ref{2n} we obtain
the following.

\begin{cor}\label{cor1}
Let $\gamma$ and $\tau$ be operator ideal norms on $\cal F$. Let
$X$ be a Banach space such that $\cal A(X)$ is $K$-amenable and
let $P:X\to X$ be a bounded projection. Set $Y=\mathrm{rg}\,P$ and
$Z=\mathrm{rg}\,(I-P)$. If $\gamma$ and $\tau$ are equivalent on
one of $\cal F(Y,Z)$ or $\cal F(Z,Y)$, say $c\gamma \leq\tau\leq
C\gamma$, then we must have $c\kappa^{-2} \gamma\leq\tau\leq
C\kappa^2 \gamma$ on one of $\cal F(Y)$ or $\cal F(Z)$, for some
$\kappa \leq 4 K^2 \|P\|\|I-P\|\max\{\|P\|^4,\|I-P\|^4\}$.
\end{cor}

\begin{proof}
By Proposition~\ref{am+ap}, at least one of the product maps
$\pi_1:\cal A(Y,Z)\widehat{\otimes}\cal A(Z,Y)\to\cal A(Z)$ or
$\pi_2:\cal A(Z,Y)\widehat{\otimes}\cal A(Y,Z)\to\cal A(Y)$ is
onto with inversion constant no greater than $4
K\|P\|\|I-P\|\max\{\|P\|^3,\|I-P\|^3\}$. To fix ideas, suppose
$\pi_1:\cal A(Y,Z)\widehat{\otimes}\cal A(Z,Y)\to\cal A(Y)$ is
onto. As $\cal A(X)$ is $K$-amenable it has a BAI of bound $K$,
and so, $\cal A(Y)$ has a BAI of bound $K\|P\|$. Now one just
needs to apply Theorem~\ref{2n}.
\end{proof}

The estimate for $\kappa$ given in the last corollary is very
unlikely to be sharp. However, to the effects of the present
paper, the significant fact about it is that it depends only on
the amenability constant and the given projection. The importance
of this fact will be fully appreciated in Section~\ref{T} when we
prove the non-amenability of $\cal A(T)$ for $T$ the
Tsirelson's~space.

Following are some important consequences of Corollary~\ref{cor1}.

In what follows, we denote by $\Gamma_p$ $(1\leq p\leq\infty)$ the
operator ideal of all bounded linear maps between Banach spaces
that factor through $\ell_p$ endowed with the operator ideal norm
\[
\gamma_p(T:X\to Y) := \inf\big\{\|R\|\|S\| :
X\stackrel{S}{\to}\ell_p\stackrel{R}{\to}Y\text{ and
}\,RS=T\big\}.
\]

Recall also that a Banach space $X$ is said to be of cotype~2 if
there exists a constant $C$ such that, for all finite subsets
$\{x_1,x_2,\ldots,x_n\}$ of $X$, we have
\[
\Big(\sum\nolimits_i \|x_i\|^2\Big)^{1/2} \leq C
2^{-n}\sum\nolimits_{t\in\{-1,1\}^n} \Big\|\sum\nolimits_i t_i
x_i\Big\|.
\]

\begin{cor}\label{cot2}
Let $X$ and $Y$ be infinite-dimensional Banach spaces with the
$\lambda$-BAP. If none of them is isomorphic to a Hilbert space
and if $X^*$ and $Y$ are both of cotype 2 then $\cal A(X\oplus Y)$
is not amenable.
\end{cor}

\begin{proof}
Since $X^*$ and $Y$ are both of cotype 2, by Pisier's abstract
version of Grothendieck's inequality \cite[Theorem~4.1]{P}, we
have that $\cal B(X,Y) = \Gamma_2(X,Y)$, and hence that
$\|\,.\,\|\sim\gamma_2$ on $\cal F(X,Y)$. Suppose towards a
contradiction that $\cal A(X\oplus Y)$ is amenable. By
Corollary~\ref{cor1}, either $\|\,.\,\| \sim \gamma_2$ on $\cal
F(X)$ or $\|\,.\,\| \sim \gamma_2$ on $\cal F(Y)$. Suppose
$\|\,.\,\| \sim \gamma_2$ on $\cal F(X)$, so, for some constant
$M$ we have $\gamma_2(T) \leq M \|T\|$ $(T\in\cal F(X))$. Let
$E\subset X$ be a finite-dimensional subspace. By \cite[\S~16.9,
Corollary]{DF}, there exists $T_E\in\cal F(X)$ such that $T_E(x) =
x$ $\,(x\in E)$ and $\|T_E\|\leq\lambda+1$. Let $\imath_E:E\to X$
be the inclusion map. Then we have
\[
\gamma_2 (I\big|_E) = \gamma_2 (T_E \imath_E) \leq \gamma_2 (T_E)
\|\imath_E\| \leq M \|T_E\| \leq M (\lambda+1).
\]
This last holds for any $E$, so, $\sup_E \gamma_2 (I\big|_E)\leq
M(\lambda+1)$. By \cite[Proposition~5.2]{LP}, $\gamma_2(I)\leq
M(\lambda+1)$, i.e., $X$ is isomorphic to a Hilbert space,
contrary to assumption.

Analogously, if $\|\,.\,\|\sim\gamma_2$ on $\cal F(Y)$, we find
that $Y$ must be isomorphic to a Hilbert space, contradicting the
hypotheses once again. Thus, $\cal A(X\oplus Y)$ cannot be
amenable.
\end{proof}

Let $1<p<\infty$. The $p$-th James space, $\frak J_p$, is the
completion of the linear space of complex sequences with finite
support in the norm
\begin{eqnarray*}
\|(\alpha_n)\|_{\frak J_p} &=&
\sup\Big\{\Big(\sum\nolimits_{n=1}^{m-1} |\alpha_{i_n} -
\alpha_{i_{n+1}}|^p\Big)^{1/p} : m,i_1,\ldots,i_m\in\Bbb N,
\\ & & \hspace{6cm} m\geq 2\text{ and }i_1<\ldots<i_m\Big\}.
\end{eqnarray*}
It is unknown if $\cal A(\frak J_p)$ is amenable for any $p$.
However, as a consequence of Corollary~\ref{cot2} we have the
following.

\begin{cor}
The algebra $\cal A(\frak J_p\oplus\frak J_p^*)$ is not amenable
for any $p\in[2,\infty[$.
\end{cor}

\begin{proof}
By \cite[Theorem~1]{P1}, $\frak J_p^*$ has cotype 2 and neither
$\frak J_p$ nor $\frak J_p^*$ is isomorphic to a Hilbert space. So
we can apply Corollary~\ref{cot2}.
\end{proof}

Recall from \cite{LP} that a Banach space $X$ is said to be an
$\cal L_p$-space if it contains a net $(X_\alpha)$ of
finite-dimensional subspaces, directed by inclusion, whose union
is dense in $X$, and such that $\sup_\alpha
d(X_\alpha,\ell_p^{\,\dim{X_\alpha}})<\infty$.

Our next result generalizes \cite[Theorem~6.9]{GJW}.

\begin{cor}\label{p,q}
Let $X$ be an $\cal L_p$-space, and let $Y$ be an $\cal
L_q$-space, where $1\leq p,q\leq\infty$. Then $\cal A(X\oplus Y)$
is amenable if and only if one of the following holds:
\begin{enumerate}
\item[$a)$] $p = q$. 

\item[$b)$] $p = 2$ and $1<q<\infty$. 

\item[$c)$] $1<p<\infty$ and $q = 2$. 
\end{enumerate}
\end{cor}

\begin{proof}
Since the direct sum of two $\cal L_p$-spaces (resp.~of an $\cal
L_p$-space with $1<p<\infty$ and an $\cal L_2$-space) is an $\cal
L_p$-space, and the algebra of approximable operators on an $\cal
L_p$-space is always amenable \cite[Theorem~6.4]{GJW}, it is clear
that if $(a)$ (resp. $(b)$ or $(c)$) is satisfied then $\cal
A(X\oplus Y)$ is amenable.

Now suppose that none of $(a)$, $(b)$ or $(c)$ is satisfied. We
want to show that $\cal A(X\oplus Y)$ is not amenable. By
\cite[Theorem~III(a)]{LR} and \cite[Corollary~5.5]{GJW}, it is
enough to consider the following two cases: (i) $p<2\leq q$, and
(ii) $p<q<2$. The case (i) follows from Corollary~\ref{cot2} above
since for $p<2$ (resp.~$2\leq q$) an $\cal L_p$-space (resp.~the
dual of an $\cal L_q$-space) has cotype 2. In dealing with the
second case we use the following result from \cite{K}, that we
state as in \cite[\S~26.5. Corollary~2]{DF}: \vspace{2mm}

\noindent{\bf Theorem} (Kwapie\'{n}). Let $1\leq p\leq r\leq
q\leq\infty$. Then $\cal B(\ell_q,\ell_p) =
\Gamma_r(\ell_q,\ell_p)$. \vspace{2mm}

Let $p<q<2$ and let $r\in\;]p,q[$. Using Kwapie\'{n}'s theorem and
\cite[Theorem~III(c)]{LR} it can be shown that there exists a
constant $M$ so that
\[
\sup\{\gamma_r(T|_E) : E\subset Y\;\text{a finite-dimensional
subspace}\} \leq M \|T\|\qquad (T\in\cal B(Y,X)).
\]
By \cite[Corollary~8.9]{P}, there is an $L_r$-space so that $\cal
B(Y,X) = \Gamma_{L_r}(Y,X)$, where $\Gamma_{L_r}(Y,X)$ denotes the
space of all operators from $X$ to $Y$ that factor through $L_r$
with the norm $\gamma_{L_r}(T) := \inf\{\|R\|\|S\| :
X\stackrel{S}{\to} L_r\stackrel{R}{\to} Y\text{ and }RS=T\}$.
Assume towards a contradiction that $\cal A(X\oplus Y)$ is
amenable. Then, by Corollary~\ref{cor1}, either $\gamma_{L_r} \sim
\|\,.\,\|$ on $\cal F(X)$ or $\gamma_{L_r} \sim \|\,.\,\|$ on
$\cal F(Y)$. We show that none of these can happen. Indeed,
suppose to fix ideas, that $\gamma_{L_r} \sim \|\,.\,\|$ on $\cal
F(X)$. Then, by \cite[Theorem~4.3]{LR}, $\gamma_{L_r}(I_X)
<\infty$, and so, $X$ is isomorphic to a complemented subspace of
an $L_r$-space which, by \cite[Theorem~III(b)]{LR}, must be an
$\cal L_r$-space. But this is impossible since $p\ne r$.
Analogously, if $\gamma_{L_r} \sim \|\,.\,\|$ on $\cal F(Y)$, we
find that $Y$ is an $\cal L_r$-space as well as an $\cal
L_q$-space reaching again the same absurd. Thus, neither
$\gamma_{L_r} \sim \|\,.\,\|$ on $\cal F(X)$ nor $\gamma_{L_r}
\sim \|\,.\,\|$ on $\cal F(Y)$. It follows that $\cal A(X\oplus
Y)$ cannot be amenable and this concludes the proof.
\end{proof}

\begin{rem}
It should be noted that the argument of \cite[Theorem~6.9]{GJW}
can be extended without difficulty to cover the more general
situation of Corollary~\ref{p,q} when $1<p,q<\infty$.
\end{rem}

We now turn our attention to the second fact mentioned at the
beginning of this section, namely, the equivalence between
surjectivity of $\cal A(\ell_p,X)\widehat{\otimes}\cal
A(X,\ell_p)\to\cal A(X)$ and openness of $\cal A(\ell_p,X)\times\cal
A(X,\ell_p)\to\cal A(X)$. It is not hard to see that the reason why
this last holds is that $\ell_p\cong \ell_p(\ell_p)$ $(1\leq
p\leq\infty)$, or more precisely, because $\gamma_p$, being a norm,
must satisfy the triangle inequality. In what remains of this
section we look at this in more detail.

Let $Z$ be an infinite-dimensional Banach space. Given any pair of
Banach spaces $(X,Y)$ we let
\[
\gamma_Z (T) := \inf\big\{\|R\|\|S\| : RS = T,\, S\in\cal
F(X,Z)\;\text{and}\;R\in\cal F(Z,Y)\big\} \quad (T\in\cal F(X,Y)).
\]

In general $\gamma_Z$ need not be a norm on $\cal F(X,Y)$. For
example, let $Z_n = \ell_p^n\oplus\ell_2$ for some
$p\in(2,\infty)$ fixed, let $I:\ell_p^{2n}\to \ell_p^{2n}$ be the
identity map, and let $P_1$ (resp.~$P_2$) be the natural
projection onto the first (resp.~last) $n$ coordinates. Then
$\gamma_{Z_n}(P_1 + P_2)$ tends to $\infty$ with $n$ while
$\gamma_{Z_n}(P_1) + \gamma_{Z_n}(P_2) = 2$ for all $n$. Indeed,
suppose towards a contradiction that $\gamma_{Z_n}(P_1 + P_2)<C$
for some constant $C$ independent of $n$. Then for every $n$ there
is $E_n\subset\ell_p^n\oplus\ell_2$ and a linear isomorphism
$T_n:\ell_p^{2n}\to E_n$ such that $\|T_n\|\|T_n^{-1}\|<C$. Let
$Q_n$ be the natural projection from $\ell_p^n\oplus\ell_2$ onto
$\ell_p^n$, and let $(x_{n,i})$ be a basis for $E_n$. Without loss
of generality, let $Q_n(x_{n,1}),\ldots,Q_n(x_{n,m})$ be a maximal
subset of linearly independent vectors from $\big\{Q_n(x_{n,i}) :
1\leq i\leq 2n\big\}$, so $m\leq n$. Taking linear combinations if
necessary, we can pass to a new basis of $E_n$,
$x_{n,1},\ldots,x_{n,m},y_{n,m+1},\ldots,y_{n,2n}$ say, in which
each $y_{n,i}\in\ell_2$. Thus, $E_n$ contains an isometric copy of
$\ell_2^n$ and in turn $\ell_p^{2n}$ contains a $C$-isomorphic
copy of $\ell_2^n$. But this last should hold for every $n$, which
is impossible by~\cite[Example~3.1]{FLM}. Thus, for big enough
$n$, $\gamma_{Z_n}$ is not a norm.

Let us say that the Banach space $Z$ has the
\emph{factorization-norm property} if for every pair of Banach
spaces, $(X,Y)$, $\gamma_Z$ is a norm on $\cal F(X,Y)$. It is
easily verified that if $Z$ has the factorization-norm property
then $\gamma_Z$ is an operator ideal norm on $\cal F$. Also note
from the example of the previous paragraph that the
factorization-norm property is an isometric~property.

\begin{cor}\label{cfr}
Let $X$ be an infinite dimensional Banach space such that the
algebra $\cal A(X)$ is $K$-amenable. Let $P$ be a bounded
projection on $X$. Set $Y = \mathrm{rg}\,P$ and $Z =
\mathrm{rg}\,(I-P)$. If both, $Y$ and $Z$, have the
fac\-tor\-iza\-tion-norm property, then at least one of the maps
\[
\varphi_Y:\cal F(Z,Y)\times\cal F(Y,Z)\to\cal F(Y) \quad\text{ or
}\quad \varphi_Z:\cal F(Y,Z)\times\cal F(Z,Y)\to\cal F(Z),
\]
is $M$-open for some constant $M\leq\kappa^2 K^2\|P\|\|I-P\|$ with
$\kappa$ as in Corollary~\ref{cor1}.
\end{cor}

\begin{proof}
Since $\cal A(X)$ is $K$-amenable it has a BAI of bound $K$. In
turn $\cal A(Y)$ has a BAI of bound $C = K\|P\|$ and $\cal A(Z)$
has a BAI of bound $c^{-1} = K\|I-P\|$. As $\cal A(Y,Z)$ is an
essential Banach left $\cal A(Z)$-module and an essential Banach
right $\cal A(Y)$-module we have, by \cite[Theorem~2.9.24]{Da},
that $\gamma_Y\leq C\|\cdot\|$ and $\gamma_Z\leq c^{-1}\|\cdot\|$
on $\cal F(Y,Z)$, so $c\gamma_Z \leq\gamma_Y\leq C\gamma_Z$
on~$\cal F(Y,Z)$. Thus, by Corollary~\ref{cor1}, we should have
$\kappa^{-2} c\gamma_Z \leq\gamma_Y\leq \kappa^2 C\gamma_Z$ on at
least one of $\cal F(Y)$ or $\cal F(Z)$ for some constant
$\kappa$. To fix ideas, suppose $\kappa^{-2} c\gamma_Z
\leq\gamma_Y\leq \kappa^2 C\gamma_Z$ holds on $\cal F(Y)$. Since
$\cal A(Y)$ is an essential $\cal A(Y)$-module we have, once again
by \cite[Theorem~2.9.24]{Da}, that $\gamma_Y\leq C\|\cdot\|$ on
$\cal F(Y)$. This last combined with $\kappa^{-2} c\gamma_Z
\leq\gamma_Y$ gives that $\gamma_Z\leq M\|\cdot\|$ on $\cal F(Y)$,
where $M\leq\kappa^2 C/c$, as desired. The case where $\kappa^{-2}
c\gamma_Z \leq\gamma_Y\leq \kappa^2 C\gamma_Z$ holds on $\cal
F(Z)$ is treated analogously.
\end{proof}

It seems difficult, in general, to determine whether or not a
given Banach space has the factorization-norm property. It is well
known, for instance, that any Banach space $Z$ such that $Z\cong
\ell_p(Z)$, in particular, any Banach space of the form
$\ell_p(E)$, where $E$ is some Banach space and $1\leq p\leq
\infty$, has the factorization-norm property (see
\cite[Proposition~1]{Jo}).

The following proposition is analogous to
\cite[Proposition~1]{Jo}. It gives a sufficient condition for a
Banach space to have the factorization-norm property.

\begin{prop}\label{fn}
Let $Z$ be an infinite dimensional Banach space such that for
every finite-dimensional subspace $E$ of $Z$ and every
$\varepsilon > 0$ there exist finite-dimensional subspaces $F$ and
$G$ of $Z$ such that $E\subset F$, $G$ is $(1 +
\varepsilon)$-complemented in $Z$, and for some set of vectors
$\{u_1,u_2\}$ forming a 1-unconditional basis of their  $\Bbb
R$-linear span we have that $d(G,F\oplus_u F) \leq 1 +
\varepsilon$, where $F\oplus_u F$ denotes the direct sum of two
copies of $F$ endowed with the norm $\|(x,y)\| :=
\big\|\|x\|_Z\,u_1 + \|y\|_Z\,u_2\big\|$ $((x,y)\in F\oplus F)$.
Then $Z$ has the factorization-norm property.
\end{prop}

\begin{proof}
Of course, only the triangle inequality needs to be verified. For
this, let $(X,Y)$ be a pair of Banach spaces, let $T_1, T_2\in
\cal F(X,Y)$, and let $\varepsilon > 0$ be arbitrary. Let
$S_i\in\cal F(X,Z)$ and $R_i\in\cal F(Z,Y)$ be such that $R_i S_i
= T_i$ and $\|R_i\|\|S_i\| \leq \gamma_Z(T_i) + \varepsilon$ $(i =
1,2)$. Set $E = S_1 X + S_2 X \subset Z$ and let $F$ and $G$ be
finite-dimensional subspaces of $Z$ as in the hypotheses. Let
$L:G\to F\oplus_u F$ be a linear isomorphism such that
$\|L\|\|L^{-1}\| \leq 1 + \varepsilon$, and let $P_G:Z\to G$ be a
projection onto $G$ of norm $\leq 1 + \varepsilon$. Let
$\overline{S}:X\to G$, $x\mapsto L^{-1}(S_1 x, S_2 x)$, and let
$\overline{R}:Z\to Y$, $z\mapsto (R_1 P_1 + R_2 P_2) L P_G z$,
where $P_1$ and $P_2$ denote the canonical coordinate projections
onto the first and second components of $F\oplus_u F$,
respectively.

It is easily seen that $\,T_1 + T_2 =
\overline{R}\,\overline{S}\,$, that
\[
\|\overline{S}\| \leq \|L^{-1}\| \big\|\|S_1\| u_1 + \|S_2\|
u_2\big\|,
\]
and that
\[
\|\overline{R}\| \leq (1 + \varepsilon) \|L\| \big\|\|R_1\| u_1^*
+ \|R_2\| u_2^*\big\|,
\]
where $u_1^*,u_2^*$ is the basis dual to $u_1,u_2$. By
\cite[Main~Lemma]{JR}, we can assume that
\[
\big\|\|S_1\| u_1 + \|S_2\| u_2\big\| \big\|\|R_1\| u_1^* +
\|R_2\| u_2^*\big\| = \|R_1\|\|S_1\| + \|R_2\|\|S_2\|.
\]
Then
\begin{eqnarray*}
\gamma_Z(T_1 + T_2) &\leq& \|\overline{R}\|\|\overline{S}\|
\;\leq\; (1 + \varepsilon)^2
(\|R_1\|\|S_1\| + \|R_2\|\|S_2\|) \\
&\leq& (1 + \varepsilon)^2 (\gamma_Z(T_1) + \gamma_Z(T_2) + 2
\varepsilon).
\end{eqnarray*}
Since $\varepsilon$ is arbitrary the desired conclusion follows.
\end{proof}

All the following Banach spaces are easily seen to satisfy the
condition of Proposition~\ref{fn} and hence have the
factorization-norm property. \vspace{1mm}

\begin{ex}
Any Banach space $Z$ with a 1-unconditional Schauder basis
$\mathbf{z} = (z_n)$ such that (i) $\,\liminf_n
d\big([z_i]_1^n,[z_i]_{n+1}^{2n}\big) = 1$; and (ii)
$\|\sum_{i=1}^{2n} a_i z_i\| = \|\sum_{i=1}^{2n} b_i z_i\|$ for all
scalar sequences $a_1, a_2, \ldots, a_{2n}$ and $b_1, b_2, \ldots,
b_{2n}$ such that $\|\sum_{i=1}^n a_i z_i\| = \|\sum_{i=1}^n b_i
z_i\|$ and $\|\sum_{i=n+1}^{2n} a_i z_i\| = \|\sum_{i=n+1}^{2n} b_i
z_i\|$ $\,(n\in\Bbb N)$. 
\end{ex}

\begin{ex}
Any Banach space of the form $(\bigoplus_k E_k)_{\mathbf{z}}$,
where $\mathbf{z}$ is as in the previous example and $E_k = E$
$(k\in\mathbb{N})$ for some Banach space $E$. 
\end{ex}

\begin{ex}
Any \emph{Johnson space} in the sense of \cite[Definition~3.1]{B}.
\end{ex}

It is unclear, however, whether or not the factorization-norm
property is an essential hypothesis in Corollary~\ref{cfr}. In fact,
the following argument suggests that the same conclusion or at least
a similar one might hold without this assumption.

Let $1<p\leq\infty$ and let $X$ be an infinite dimensional Banach
space such that $\cal A(X)$ is $K$-amenable. Then $\cal
A(\ell_p(X)) (= \cal A_0(\ell_p(X)))$ is $K$-amenable as well. Let
$P:X\to X$ be a bounded projection. Set $Y = \text{rg}\,P$ and $Z
= \text{rg}\,(I-P)$. As $\ell_p(Y)$ and $\ell_p(Z)$ have the
factorization-norm property and $\|\bigoplus_{i=1}^\infty
P:\ell_p(X)\to\ell_p(Y)\| = \|P\|$ there is, by
Corollary~\ref{cfr}, a constant $M$, depending only on $K$ and
$P$, so that, at least one of the maps,
\[
\varphi_{\ell_p(Y)}:\cal F(\ell_p(Z),\ell_p(Y))\times\cal
F(\ell_p(Y),\ell_p(Z))\to\cal F(\ell_p(Y)),
\]
or
\[
\varphi_{\ell_p(Z)}:\cal F(\ell_p(Y),\ell_p(Z))\times\cal
F(\ell_p(Z),\ell_p(Y))\to\cal F(\ell_p(Z)),
\]
is $M$-open. To fix ideas, suppose $\varphi_{\ell_p(Y)}$ is
$M$-open. Then
\[
\varphi_Y:\cal F(\ell_p(Z),Y)\times\cal F(Y,\ell_p(Z))\to\cal
F(Y),
\]
is $M$-open too.

\begin{rem}
The fact that at least one of the maps $\varphi_{\ell_p(Y)}$ or
$\varphi_{\ell_p(Z)}$ above is $M$-open if $\cal A(X)$ is
amenable, still holds for $p=1$, but this case needs to be treated
separately as $\cal A(\ell_1(X))\ne\cal A_0(\ell_1(X))$ (see
Lemma~\ref{p=1} below).
\end{rem}

Now let $(T_\alpha)$ be a BAI for $\cal F(Y)$. Then, by the above,
we have that for every $1\leq p\leq\infty$ and every $\alpha$
there are operators $\,S_{\alpha,p}:Y\to \ell_p(Z)\,$ and
$\,R_{\alpha,p}:\ell_p(Z)\to Y\,$ such that $R_{\alpha,p}\,
S_{\alpha,p} = T_\alpha$ and $\|R_{\alpha,p}\|\|S_{\alpha,p}\|
\leq M \|T_\alpha\|$. The fact that $M$ is independent of $\alpha$
and $p$ suggests the following might be~true: \vspace{3mm}

\noindent $\blacksquare\,$ There exist $k\in\Bbb N$ and a positive
constant $\widetilde{M}$ such that for every index $\alpha$, there
are operators $\,R_\alpha:\bigoplus_{i=1}^k Z\to Y\,$ and
$\,S_\alpha:Y\to\bigoplus_{i=1}^k Z\,$ so that $R_\alpha S_\alpha
= T_\alpha$ and $\|R_\alpha\|\|S_\alpha\|\leq
\widetilde{M}\|T_\alpha\|$, that is, the product map
\[
\varphi:\cal F\Big(\bigoplus\nolimits_{i=1}^k Z,Y\Big)\times\cal
F\Big(Y,\bigoplus\nolimits_{i=1}^k Z\Big)\to \cal F(Y)
\]
is $\widetilde{M}$-open.

\section{Tsirelson-like spaces}\label{T}

As announced earlier, in the final section of this paper we
establish the non-amenability of the algebra of approximable
operators on the Tsirelson space. In fact, we shall obtain this as a
consequence of a more general result (see Theorem~\ref{CO} below).

We start with a definition. It is closely related to the old
notion of crude finite representability introduced in \cite{Ja}.

\begin{definition}\label{dcfr}
A Banach space $Y$ is said to be \emph{$M$-crudely $\pi$-finitely
representable} in a Banach space $Z$ if for every
finite-dimensional subspace $E$ of $Y$, there exist a finite-rank
projection $P:Y\to Y$ whose range contains $E$, and operators
$S:Y\to Z$ and $R:Z\to Y$ such that $R S = P$ and $\|R\|\|S\|\leq
M$.
\end{definition}

This last definition is justified by the following.

\begin{lem}\label{ef}
Let $Y$ be a $\pi_1$-space and let $Z$ be a Banach space. Then all
the following are equivalent: \vspace{1mm}
\begin{itemize}
\item[$a)$] $Y$ is $(M+\varepsilon)$-crudely $\pi$-finitely
representable in $Z$ for every $\varepsilon>0$. \vspace{1mm}

\item[$b)$] $\cal F(Z,Y)\times\cal F(Y,Z) \to\cal F(Y)$ is
$(M+\varepsilon)$-open for every $\varepsilon>0$. \vspace{1mm}

\item[$c)$] For every $\varepsilon>0$ there exists in $\cal F(Y)$
a bounded net of projections, $(P_\alpha)$, converging strongly to
the identity operator on $Y$, and such that $\sup_\alpha \gamma_Z
(P_\alpha)\leq M+\varepsilon$, that is, such that $Z$ contains
$P_\alpha Y$'s uniformly $(M+\varepsilon)$-complemented.
\end{itemize}
\end{lem}

\begin{proof}
It is easy to see that $(a)$ $\Rightarrow$ $(b)$ and $(b)$
$\Rightarrow$ $(c)$. That $(c)$ $\Rightarrow$ $(a)$ follows from
\cite[Lemma~2.4]{JRZ}.
\end{proof}

\begin{rem}
If we simply assume in the last lemma that $Y$ is a
$\pi_\lambda$-space, then still $(a)$ $\iff$ $(c)$ and $(a)$
$\Rightarrow$ $(b)$, but $(b)$ implies that $Y$ is $(\lambda M +
\varepsilon)$-crudely $\pi$-finitely representable in $Z$ for
every $\varepsilon>0$.
\end{rem}

Before passing to the main result of the section we need another
result that we collect as a lemma.

Let $E$, $F$ be Banach spaces and let $1\leq p\leq\infty$. Recall
that we have identified $\ell_p^m$ with the linear span of the
first $m$ vectors of the unit vector basis of $\ell_p$, so we have
a natural linear isometry $\cal
A(\ell_p^m(E),\ell_p^k(F))\hookrightarrow \cal
A(\ell_p^n(E),\ell_p^l(F))$, whenever $n\geq m$ and $l\geq k$. We
denote by $\cal A_0(\ell_p(E),\ell_p^k(F))$ (resp.~$\cal
A_0(\ell_p^m(E),\ell_p(F))$) the inductive limit of the direct
system formed by the spaces $\cal A(\ell_p^m(E),\ell_p^k(F))$
$(m\in\Bbb N)$ (resp.~$\cal A(\ell_p^m(E),\ell_p^k(F))$ $(k\in\Bbb
N)$) together with the corresponding isometric embeddings. There
are also natural linear isometries $\cal
A_0(\ell_p(E),\ell_p^k(F))\hookrightarrow \cal
A_0(\ell_p(E),\ell_p^l(F))$ and $\cal
A_0(\ell_p^m(E),\ell_p(F))\hookrightarrow \cal
A_0(\ell_p^n(E),\ell_p(F))$ $(l\geq k,\,n\geq m)$. We denote by
$\cal A_0(\ell_p(E),\ell_p(F))$ the common inductive limit of the
direct systems formed by $\{\cal A_0(\ell_p(E),\ell_p^k(F)) :
k\in\Bbb N\}$ and $\{\cal A_0(\ell_p^m(E),\ell_p(F)) : m\in\Bbb
N\}$ and their respective families of isometric embeddings. It is
not hard to see that $\cal A_0(\ell_p(E),\ell_p(F))$ is also the
inductive limit of the direct system formed by all spaces $\cal
A_0(\ell_p^m(E),\ell_p^k(F))$ and the isometric embeddings $\cal
A(\ell_p^m(E),\ell_p^k(F))\hookrightarrow \cal
A(\ell_p^n(E),\ell_p^l(F))$ $(n\geq m,\,l\geq k)$.

\begin{lem}\label{p=1}
Let $1\leq p\leq\infty$ and let $X$ be a Banach space such that
$\cal A(X)$ is $K$-amenable.  Let $P:X\to X$ be a bounded
projection. Set $Y=\mathrm{rg}\,P$ and $Z=\mathrm{rg}\,(I-P)$.
Then at least one of the maps
\[
\varphi_1:\cal A_0(\ell_p(Z),\ell_p(Y))\times\cal
A_0(\ell_p(Y),\ell_p(Z))\to\cal A_0(\ell_p(Y)),
\]
or
\[
\varphi_2:\cal A_0(\ell_p(Y),\ell_p(Z))\times\cal
A_0(\ell_p(Z),\ell_p(Y))\to\cal A_0(\ell_p(Z)),
\]
is $(M+\delta)$-open for every $\delta>0$ and some constant $M$
that depends only on $K$ and $\|P\|$.
\end{lem}

\begin{proof}
Let $X$ be a Banach space such that $\cal A(X)$ is $K$-amenable.
Then $\cal A(\ell_p^n(X))$ is $K$-amenable for every $n\in\Bbb N$
and every $1\leq p\leq\infty$. Let $P$, $Y$ and $Z$ be as in the
hypotheses. As $\|\bigoplus_{1}^n P:\ell_p^n(X)\to\ell_p^n(Y)\| =
\|P\|$ $(n\in\Bbb N)$, there exists, by Proposition~\ref{am+ap}, a
constant $M = M(K,\|P\|)$ so that for each $n\in\Bbb N$ at least
one of the product maps
\begin{equation}\label{5}
\pi_{1,n}:\cal A(\ell_p^n(Z),\ell_p^n(Y))\widehat{\otimes}\cal
A(\ell_p^n(Y),\ell_p^n(Z))\to\cal A(\ell_p^n(Y)),
\end{equation}
or
\begin{equation}\label{6}
\pi_{2,n}:\cal A(\ell_p^n(Y),\ell_p^n(Z))\widehat{\otimes}\cal
A(\ell_p^n(Z),\ell_p^n(Y))\to\cal A(\ell_p^n(Z)),
\end{equation}
is onto with inversion constant no greater than $M$.

Without loss of generality, assume there is an increasing sequence
of positive integers, $(n_k)$, so that $\pi_{1,n_k}$ is onto with
inversion constant no greater than $M$ for all $k$. Fix $k\in\Bbb
N$, let $\varepsilon>0$ and let $T\in\cal A(\ell_p^{n_k}(Y))$.
There is $\sum_i R_i\otimes S_i\in \cal
A(\ell_p^{n_k}(Z),\ell_p^{n_k}(Y))\widehat{\otimes}\cal
A(\ell_p^{n_k}(Y),\ell_p^{n_k}(Z))$ such that $\sum_i R_i S_i = T$
and $\sum_i \|R_i\|\|S_i\|\leq(M+\varepsilon)\|T\|$. Moreover, we
can assume $\lim_i \|R_i\| = 0 = \lim_i \|S_i\|$. For each
$i\in\Bbb N$, let $P_i$ denote the coordinate projection of
$\ell_p(\ell_p^{n_k}(Y))$ onto its $i$-th summand, and let
$\imath_i$ denote the embedding of the $i$-th summand
into~$\ell_p(\ell_p^{n_k}(Y))$). Let $R = \sum_i R_i P_i$ and $S =
\sum_i \imath_i S_i$. Then $RS = T$, $\|R\|^q \leq \sum_i
\|R_i\|^q$ (resp.~$\leq\max_i \|R_i\|$ if $q=\infty$) and $\|S\|^p
\leq \sum_i \|S_i\|^p$ (resp.~$\leq\max_i \|S_i\|$ if $p=\infty$).
Furthermore, by \cite[\S 2]{JR}, we can choose the $R_i$'s and
$S_i$'s in such a way that $\|R\|\|S\|\leq(1+\varepsilon) \sum_i
\|R_i\|\|S_i\|$, and hence $\|R\|\|S\|\leq
(1+\varepsilon)(M+\varepsilon)\|T\|$. As this last holds for
arbitrary $\varepsilon$, the bilinear map
\[
\varphi_{n_k}:\cal A_0(\ell_p(Z),\ell_p^{n_k}(Y))\times\cal
A_0(\ell_p^{n_k}(Y),\ell_p(Z))\to\cal A(\ell_p^{n_k}(Y))
\qquad(k\in\Bbb N),
\]
is $(M+\delta)$-open for any $\delta>0$.

Now let $\cal T\in\cal A_0(\ell_p(Y))$ and $\varepsilon>0$. There
exists a sequence $(T_k)$ in $\cal A_0(\ell_p(Y))$ such that
$T_k\in\cal A(\ell_p^{n_k}(Y))$ for every $k$, $\sum_k T_k = \cal
T$ and $\sum_k \|T_k\|\leq \|\cal T\|+\varepsilon$. By the
previous part, there exist $R_k\in\cal
A_0(\ell_p(Z),\ell_p^{n_k}(Y))$ and $S_k\in\cal
A_0(\ell_p^{n_k}(Y),\ell_p(Z))$ such that $R_k S_k = T_k$ and
$\|R_k\|\|S_k\|\leq (M+\varepsilon)\|T_k\|$ $(k\in\Bbb N)$. Let
$\pi_k$ denote the projection of $\ell_p(\ell_p(Z))$ onto its
$k$-th summand, and let $\jmath_k$ denote the natural embedding of
the $k$-th summand back into~$\ell_p(\ell_p(Z))$. Define $\cal R =
\sum_k R_k \pi_k$ and $\cal S = \sum_k \jmath_k S_k$. Then $\cal R
\cal S = \cal T$ and an argument similar to the one of the
previous paragraph shows that $\|\cal R\|\|\cal S\|\leq
(1+\varepsilon)(M+\varepsilon)(\|\cal T\|+\varepsilon)$. As
$\ell_p(\ell_p(Z))\cong \ell_p(Z)$, it follows that the product
\[
\varphi_1:\cal A_0(\ell_p(Z),\ell_p(Y))\times\cal
A_0(\ell_p(Y),\ell_p(Z))\to\cal A_0(\ell_p(Y)),
\]
is $(M+\delta)$-open for any $\delta>0$.
\end{proof}

\begin{rem}
When $p>1$ things are much simpler. Indeed, in these cases the
claim of the lemma can be easily obtained from
Corollary~\ref{cfr}.
\end{rem}

\begin{thm}\label{CO}
Let $X$ be a Banach space with a 1-unconditional basis $(x_i)$.
Suppose there exist $\delta\geq 1$ and $1\leq p\leq\infty$ such
that \vspace{1mm}
\begin{itemize}
\item[(i)] For each $n\in\mathbb{N}$ there is $m\in\mathbb{N}$ so
that if $\,F\subseteq [x_i]_{i=m}^\infty\,$ is a subspace spanned
by $n$ disjointly supported vectors then $\,\inf\big\{d(F,E):
E\;\text{is a subspace of $\ell_p$}\big\}\leq\delta$. \vspace{1mm}

\item[(ii)] $\,\inf\big\{d([x_i]_{i=1}^n,E) : E \text{ is a
subspace of $\ell_p$}\big\}\underset{n}{\to} \infty$. \vspace{1mm}
\end{itemize}
\noindent Then $\cal A(X)$ is not amenable.
\end{thm}

\begin{proof}
Suppose towards a contradiction that $\cal A(X)$ is $K$-amenable
for some $K\geq 1$. Let $\delta \geq 1$ and $1\leq p\leq\infty$ be
as in the hypotheses. By Lemma~\ref{p=1}, for every $m\in\Bbb N$
at least one of the maps
\[
\varphi_{1,m}:\cal A_0\big(\ell_p([x_i]_{m+1}^\infty),
\ell_p([x_i]_{1}^m)\big)\times\cal
A_0\big(\ell_p([x_i]_{1}^m),\ell_p([x_i]_{m+1}^\infty)\big)\to
\cal A_0\big(\ell_p([x_i]_{1}^m)\big),
\]
or
\[
\varphi_{2,m}:\cal A_0\big(\ell_p([x_i]_{1}^m),
\ell_p([x_i]_{m+1}^\infty)\big)\times\cal
A_0\big(\ell_p([x_i]_{m+1}^\infty),\ell_p([x_i]_{1}^m)\big)\to
\cal A_0\big(\ell_p([x_i]_{m+1}^\infty)\big),
\]
is $M$-open, where $M$ depends only on $K$ and the norm of the
natural projection $P_m:X\to [x_i]_{i=1}^m$. But $\varphi_{2,m}$
cannot be open since otherwise, by Lemma~\ref{ef},
$\,\ell_p([x_i]_{m+1}^\infty)\,$ would be crudely $\pi$-finitely
representable in $\,\ell_p([x_i]_{1}^m)\simeq\ell_p$, which is
impossible by (ii). Thus, $\varphi_{1,m}$ is $M$-open and by
Lemma~\ref{ef}, $\ell_p([x_i]_{1}^m)$ is $(M+1)$-crudely
$\pi$-finitely representable in $\ell_p([x_i]_{m+1}^\infty)$. Note
that, as $\|P_m\|=1$ for every $m$, $M$ is also independent of
$m$. Thus, for every $m\in\mathbb{N}$ and every $1\leq n\leq m$
there exist operators
\[
S_{m,n}:[x_i]_{i=1}^n\to \ell_p([x_i]_{m+1}^\infty) \quad\text{
and }\quad R_{m,n}:\ell_p([x_i]_{m+1}^\infty)\to [x_i]_{i=1}^n,
\]
so that $R_{m,n}\, S_{m,n} = I_{[x_i]_{1}^n}$ and
$\|R_{m,n}\|\|S_{m,n}\| \leq M+1$.

What remains follows closely the proof of \cite[Prop.~VI.b.3]{CS}.
Let $F_{m,n} := S_{m,n}([x_i]_{1}^n)$. By \cite[Prop.~V.6]{CS},
there are disjointly supported vectors $\,y_1, y_2, \ldots, y_N$
in $\ell_p([x_i]_{m+1}^\infty)\,$ such that $F_{m,n}$ is
$2$-isomorphic to a subspace of $\,F := [y_j]_{1}^N$. Clearly, we
can assume that all the $y_j$'s have finite support. Let $P_k$ be
the projection of $\ell_p([x_i]_{m+1}^\infty)$ onto its $k$-th
summand. By (i), there exists $m$ so that $\,\inf\big\{d([P_k
y_j]_{j=1}^N,E) : E\; \text{a subspace of $\ell_p$}\big\}
\leq\delta$ for every $k\in\Bbb N$. As
$\,F\subseteq\big(\bigoplus_k\,[P_k y_j]_{j=1}^N\big)_p$, it
follows from this last and the fact that
$\ell_p\cong\ell_p(\ell_p)$, that $\,\inf\big\{d(F,E) : E\;\text{a
subspace of $\ell_p$}\big\}\leq \delta$, and in turn that
$\,\inf\big\{d([x_i]_{1}^n,E) : E\; \text{a subspace of
$\ell_p$}\big\}\leq 2 \delta C$, contradicting (ii). Thus $\cal
A(X)$ cannot be amenable.
\end{proof}

We apply Theorem~\ref{CO} to a class of `Tsirelson-like' spaces
introduced in \cite{FJ}, which contains the dual of the original
Tsirelson's space as a particular case.

Let us recall briefly the definition of the dual of Tsirelson's
space, $T$, as given in \cite[\S~2]{FJ}. Let $(t_n)$ denote the
unit vector basis of $c_{00}$ (the space of scalar sequences with
finite support). If $E$, $F$ are finite, non-empty subsets of
$\mathbb{N}$, we write $E<F$ to mean that $\max E < \min F$. For
any $E\subset\mathbb{N}$ and any $x = \sum_n \alpha_n t_n \in
c_{00}$, define $Ex := \sum_{n\in E} \alpha_n t_n$. Set
$\|\,.\,\|_0 := \|\,.\,\|_{c_0}$ and, for $m\geq 0$, define
\[
\|x\|_{m+1} := \max\bigg\{\|x\|_m, 2^{-1}\max\Big[
\sum\nolimits_{j=1}^k \|E_j x\|_m\Big]\bigg\}\qquad (x\in c_{00}),
\]
where the inner maximum is taken over all possible choices of
finite subsets $E_1, E_2, \ldots, E_k$ of $\mathbb{N}$, such that:
$\{k\}\leq E_1 < E_2 < \ldots < E_k$. It is easily verified that
$\|\,.\,\|_m$ is a norm on $c_{00}$ for every $m$, and that, for
each $x\in c_{00}$ the sequence $(\|\,.\,\|_m)$ is non-decreasing
and majorized by $\|x\|_{\ell_1}$. Thus we can define
\[
\|x\| := \lim_{m\to\infty} \|x\|_m \qquad (x\in c_{00}).
\]
The latter is a norm on $c_{00}$. The dual $T^*$ of Tsirelson's
space is defined as the completion of $c_{00}$ in the last norm.
It is well known that the sequence $(t_n)$ is a normalized
1-unconditional basis for $T^*$.

For $1\leq p <\infty$, $T^{(p)}$ is defined as the set of all $\,x
= \sum_n \alpha_n t_n\,$ such that $\,\sum_n |\alpha_n|^p t_n\in
T^*$, endowed with the norm
\[
\|x\|_{(p)} = \bigg\|\sum_n |\alpha_n|^p t_n
\bigg\|^\frac{1}{p}\qquad (x\in T^{(p)}).
\]

When $1<p<\infty$, $T^{(p)}$ is the so called $p$-convexified
Tsirelson's space. Clearly, $T^{(1)}$ is nothing but $T^*$ itself.

Many important facts about $T^*(= T^{(1)})$ are shared by the
$p$-convexified Tsirelson's spaces. Among them we have the
following: \vspace{1mm}

$a)$ Each $T^{(p)}$ is reflexive (actually, they are all uniformly
convex for $p>1$); \vspace{1mm}

$b)$ Each $T^{(p)}$ contains $\ell_p^n$'s uniformly $(1<p<
\infty)$; \vspace{1mm}

$c)$ No $T^{(p)}$ contains an isomorphic copy of $\ell_r$ $(1\leq
r\leq\infty)$. \vspace{1mm}

Moreover, for every $p\geq 1$ the norm on $T^{(p)}$ satisfies
\begin{equation}\label{fTp}
\|x\|_{(p)} = \max\bigg\{\|x\|_0, 2^{-\frac{1}{p}}\sup
\Big[\sum\nolimits_{j=1}^k \|E_j
x\|_{(p)}^p\Big]^\frac{1}{p}\bigg\} \qquad (x\in T^{(p)}),
\end{equation}
where the inner supremum is taken over all choices of finite
subsets $E_1, E_2, \ldots, E_k$ of $\mathbb{N}$, such that:
$\{k\}\leq E_1 < E_2 < \ldots < E_k$. Property $(b)$ above follows
easily from (\ref{fTp}).

We need one more fact about these spaces that we collect in the
next lemma.

\begin{lem}\label{*}
Let $1\leq p<\infty$. Then $T^{(p)}$ is not crudely finitely
representable in $\ell_p$.
\end{lem}

As explained in \cite[VI.B]{CS}, this follows on combining results
of Lindenstrauss and Pelczynski \cite{LP}, and Lindenstrauss and
Rosenthal \cite{LR}. Precisely, it follows from \cite[Remark after
Prop.~5.2]{LP} (see \cite[Corollary~8.9]{P} for a proof of this)
that if $X$ is a Banach space complemented in its bidual such that
for some $1\leq p<\infty$, $\sup\{\gamma_p(I_E) : \text{$E$ a
finite-dimensional subspace of $X$}\} <\infty$, then $X$ is
isomorphic to a complemented subspace of an $L_p(\mu)$ space. On
the other hand, if $X$ is a complemented subspace of an $L_p(\mu)$
space, which is not isomorphic to a Hilbert space, then it must be
an $\cal L_p$ space \cite[Theorem~III(b)]{LR}, and hence it must
contain a complemented subspace isomorphic to $\ell_p$
\cite[Proposition~7.3]{LP}. Thus, if $\,T^{(p)}\,$ $(1\leq
p<\infty)$ were crudely finitely representable in $\ell_p$ then it
would embed complementably in some $L_p(\mu)$. But $T^{(p)}$
contains $\ell_p^n$'s uniformly, and so it would contain a
complemented copy of $\ell_p$, which is an absurd.

\begin{cor}\label{T}
The algebra $\cal A(T^{(p)})$ $(1\leq p<\infty)$ is not amenable.
\end{cor}

\begin{proof}
We simply note that $T^{(p)}$ satisfies conditions (i) and (ii) of
Theorem~\ref{CO} $(1\leq p<\infty)$. Indeed, (i) is an immediate
consequence of \cite[Proposition~7.3]{BCLT} and (\ref{fTp}) above,
while (ii) follows easily from Lemma~\ref{*}.
\end{proof}

\begin{rem}
By \cite[Corollary~5.5]{GJW}, $\cal A(T)$ cannot be amenable
either.
\end{rem}

\end{document}